\input ./style/arxiv-general.cfg
\documentclass[aap,MSNbibl,dvips]{arximspdf}
\makeatletter
   \@ifpackageloaded{graphicx}{}{\usepackage{graphicx}}
\makeatother
\usepackage{mathbh}

\doi{10.1214/15-AAP1129}
\volume{26}
\issue{3}
\pubyear{2016}
\firstpage{1659}
\lastpage{1697}
\docsubty{FLA}

\makeatletter

\newcommand{\eqref}[1]{(\ref{#1})}
\newtheorem{theorem}{Theorem}[section]
\newproclaim{conj}[theorem]{Conjecture}
\newtheorem{lemma}[theorem]{Lemma}
\newtheorem{proposition}[theorem]{Proposition}
\newtheorem{corollary}[theorem]{Corollary}

\newproclaim{rem}[theorem]{Remark}

\newtheorem{lemme}{Lemma}
\newproclaim{definition}{Definition}
\makeatother

\begin{document}
\begin{frontmatter}

\title{Beyond universality in random matrix theory\thanksref{T1}}
\runtitle{Beyond universality in random matrix theory}

\begin{aug}
\author[A]{\fnms{Alan}~\snm{Edelman}\thanksref{M1}\ead[label=e1]{edelman@math.mit.edu}},
\author[A]{\fnms{A.}~\snm{Guionnet}\thanksref{M1}\corref{}\ead[label=e2]{guionnet@math.mit.edu}}
\and
\author[B]{\fnms{S.}~\snm{P\'ech\'e}\thanksref{M2}\ead[label=e3]{sandrine.peche@univ-paris-diderot.fr}}
\runauthor{A. Edelman, A. Guionnet and S. P\'ech\'e}
\affiliation{Massachusetts Institute of Technology\thanksmark{M1} and
University of Paris Diderot (Paris 7)\thanksmark{M2}}
\address[A]{A. Edelman\\
A. Guionnet\\
Department of Mathematics\\
Massachusetts Institute of Technology\\
Cambridge, Massachusetts 02139-4307\\
USA\\
\printead{e1}\\
\phantom{E-mail:\ }\printead*{e2}}
\address[B]{S. P\'ech\'e\\
U.F.R. de Math\'{e}matiques\\
University of Paris Diderot (Paris 7)\\
5 Rue Thomas Mann\\
75013 Paris\\
France\\
\printead{e3}}
\end{aug}
\thankstext{T1}{Supported by NSF Grants DMS-13-12831, DMS-10-16125, DMS-10-16086,
DMS-13-07704 and of the Simons foundation.}

%
\received{\smonth{7} \syear{2014}}
%
\revised{\smonth{7} \syear{2015}}

%
\begin{abstract}
In order to have a better understanding of finite random matrices with
non-Gaussian
entries, we study the $1/N$ expansion of local eigenvalue statistics in
both the
bulk and at the hard edge of the spectrum of random matrices. This
gives valuable
information about the smallest singular value not seen in universality laws.
In particular, we show the dependence on the fourth moment (or the
kurtosis) of the entries.
This work makes use of the so-called complex Gaussian divisible
ensembles for both Wigner and sample covariance matrices.
\end{abstract}

%
\end{frontmatter}

\section{Beyond universality}

The desire to assess the applicability of universality results in
random matrix theory
has pressed the need
to go beyond universality, in particular the need to understand the
influence of finite $n$
and what happens if the matrix deviates from Gaussian normality. In
this article, we provide exact asymptotic correction
formulas for the smallest singular value of complex matrices and bulk
statistics for complex
Wigner matrices.

``Universality,'' a term encountered in statistical mechanics, is widely
found in the field of random matrix theory.
The universality principle loosely states that eigenvalue statistics of
interest will
behave asymptotically as if the matrix elements were Gaussian.
The spirit of the term is that the eigenvalue statistics will not care
about the details of the matrix elements.

It is important to extend our knowledge of random matrices beyond
universality. In particular,
we should understand the role played by:
\begin{itemize}
\item finite $n$ and
\item non-Gaussian random variables.
\end{itemize}
From an application viewpoint, it is very valuable to have an estimate
for the departure
from universality. Real problems require that $n$ be finite, not
infinite, and it has long
been observed computationally that $\infty$ comes very fast in random
matrix theory.
The applications beg to know how fast. From a theoretical viewpoint,
there is much to
be gained in searching for proofs that closely follow the underlying
mechanisms of the mathematics.
We might distinguish ``mechanism oblivious'' proofs whose bounds require
$n$ to be well outside imaginably useful ranges, with ``mechanism
aware'' proofs that
hold close to the underlying workings of random matrices. We encourage such
``mechanism aware'' proofs.

In this article, we study the influence of the fourth cumulant on the
local statistics of the eigenvalues of random matrices of Wigner and
Wishart type.

On one hand, we study the asymptotic expansion of the smallest
eigenvalue density of large random sample covariance matrices. The
behavior of smallest eigenvalues of sample covariance matrices when
$p/n $ is close to one (and more generally) is somewhat well understood
now. We refer the reader to \cite{Edelman,Tw,Forrester,BaiYin,BAP}.
The impact of the fourth cumulant of the entries is of interest here;
we show its contribution to the distribution function of the smallest
eigenvalue density of large random sample covariance matrices
as an additional error term of order of the inverse of the dimension
(see Theorem~\ref{maintheo2}).

On the other hand, we consider the influence of the fourth moment in
the local fluctuations in the bulk. Here, we consider
Wigner matrices and discuss a conjecture of Tao and Vu \cite{TaoVu}
that the fourth moment brings a correction to the fluctuation of the
expectation of the eigenvalues in the bulk of order of the inverse of
the dimension. We prove (cf. Theorem~\ref{theotaovu}) that the
quantiles of the one point correlation function fluctuate according to
the formula predicted by Tao and Vu for the fluctuations of the
expectation of the eigenvalues.

In both cases, we consider the simplest random matrix ensembles that
are called Gaussian divisible, that is whose entries can be described
as the convolution of a distribution by the Gaussian law.
To be more precise, we consider the so-called Gaussian-divisible
ensembles, also known as Johansson--Laguerre and Johansson--Wigner
ensembles. These ensembles, defined hereafter, have been first
considered in \cite{Joh} and have the remarkable property that the
induced joint eigenvalue density can be computed. It is given in terms
of the
Itzykson--Zuber--Harich--Chandra integral. From such a formula, saddle
point analysis allows to study the local statistics of the eigenvalues.
In \cite{GTT}, this idea was used to bound the rate of convergence of
the partition function of Gaussian divisible ensemble toward their
limit by the inverse of the dimension to the power $2/3$. We precise
this study by showing that at the hard edge or in the bulk, this error
is in fact of the order of the inverse of the dimension and give the
explicit form of this error, and in particular its dependency on the
fourth moment.
It turns out that in both cases under study, the
contribution of the fourth moment to the local statistics can be
inferred from the fluctuations of the one-point correlation function,
that is of the mean linear statistics of Wigner and Wishart random
matrices. The covariance of the latter is well known, since \cite{KKP},
to depend on the fourth moments, from which our results follow.

\section{Discussion and simulations}

\subsection{Preliminaries: Real kurtosis}

We will only consider distributions whose real and imaginary parts are
independent and are identically distributed.

\begin{definition}\label{defkurt}
The \textit{kurtosis} of a distribution is
\[
\gamma=\frac{\kappa_4^{\Re}}{\sigma_{\Re}^4}=\frac{\mu
_4}{\sigma_\Re^4}-3,
\]
where $\kappa_4^{\Re}$ is the fourth cumulant of the real part,
$\sigma
_{\Re}^2$ is the variance of the real part, and $\mu_4$ is the fourth
moment about the mean.
The fourth cumulant of a centered complex distribution $P$ with i.i.d.
real and complex part with variance $\sigma_\Re^2$, is given by
\[
\kappa_4=\int\bigl|zz^*\bigr|^2 \,dP(z)-8\sigma_\Re^4=2
\kappa_4^\Re=2\gamma \sigma _{\Re}^4.
\]
\end{definition}

\subsection*{Note}
From a software viewpoint, commands such as \verb+randn+ make it
natural to take the real and the imaginary
parts to separately have mean $0$, variance $1$, and also to consider
the real kurtosis.

Example of Kurtoses $\gamma$ for distributions with mean $0$, and
$\sigma^2=1$ is provided in Table~\ref{tabb}.
\begin{table}[b]
\textwidth=223pt
\caption{Standard Kurtoses and codes}\label{tabb}
\begin{tabular*}{223pt}{@{\extracolsep{\fill}}lcc@{}}
\hline
\textbf{Distribution} & $\bolds{\gamma}$ & \textbf{Univariate code} \\
\hline
Normal & \phantom{$-$}0\phantom{.2} & {\tt randn} \\
Uniform $[-\sqrt{3},\sqrt{3}]$ & $-1.2$ &
{\tt (rand-0.5)*sqrt(12)} \\
Bernoulli & $-2$\phantom{.2} &
{\tt sign(randn)} \\
Gamma & \phantom{$-$}6\phantom{.2} & {\tt rand(Gamma()) - 1} \\
\hline
\end{tabular*}
\end{table}

For the matrices themselves, we compute the smallest
eigenvalues of the Gram matrix constructed
from $(n+\nu) \times n$ complex random matrices with
Julia \cite{julia} code provided for the reader's convenience in
Table~\ref{tabb2}.

\begin{table}[t]
\caption{Matrix codes}\label{tabb2}
\begin{tabular*}{\textwidth}{@{\extracolsep{\fill}}lc@{$\!$}}
\hline
\textbf{RM} & \textbf{Complex matrix code} \\
\hline
Normal & {\tt randn(n+$\nu$,n)+im*randn(n+$\nu$,n)} \\
Uniform & {\tt((rand(n+$\nu$,n)-0.5)+im*rand(n+$\nu$,n)-0.5))*sqrt(12)}\\
Bernoulli &
{\tt sign(randn(n+$\nu$,n))+im*sign(randn(n+$\nu$,n))} \\
Gamma &
{\tt(rand(Gamma(),n+$\nu$,n)-1)+im*(rand(Gamma(),n+$\nu$,n)-1)}\\
\hline
\end{tabular*}
\end{table}

\subsection{Smallest singular value experiments}

Let $A$ be a random $n+\nu$ by $n$ complex matrix with i.i.d. real and
complex entries all with
mean $0$, variance $1$ and kurtosis $\gamma$.
In the next several subsections, we display special cases of our
results, with experiment vs. theory curves
for $\nu=0,1$ and $2$.

We consider the cumulative distribution function
\[
F(x)=\mathbb{P} \biggl( \lambda_{\min}\bigl(AA^*\bigr) \le
\frac{x}{n} \biggr)= \mathbb{P} \biggl( \bigl( \sigma_{\min} (A)
\bigr)^2\le\frac{x}{n} \biggr),
\]
where $\sigma_{\min}(A)$ is the smallest singular value of $A$.
We also consider the density
\[
f(x)=\frac{d}{dx} F(x).
\]


In the plots to follow, we took a number of cases when $n=20,40$ and
sometimes $n=80$. We computed 2{,}000{,}000 random
samples on each of 60 processors using Julia \cite{julia}, for a total
of 120{,}000{,}000 samples of each experiment.
The runs used 75\% of the processors on a machine equipped with 8 Intel
E7-8850-2.0 GHz-24M-10 Core Xeon MP Processors.
This scale experiment, which is made easy by the Julia system, allows
us to obtain visibility on the higher order terms that would be hard to
see otherwise.
Typical runs took about an hour for $n=20$, three hours for $n=40$ and
twelve hours for $n=80$.

We remark that we are only aware of two or three instances where
parallel computing has been used in random matrix experiments.
Working with Julia is pioneering in showing just how easy this can be,
giving the random matrix experimenter a new tool
for honing in on phenomena that would have been nearly impossible to
detect using conventional methods.

\subsection{Example: Square complex matrices (\texorpdfstring{$\nu=0$}{$nu=0$})}\label{sec23}
Consider taking, a $20$ by $20$ random matrix with independent real and
imaginary entries that are uniformly distributed on $[-\sqrt{3},\sqrt{3}]$.
\[
\tt((rand(20,20)-0.5)+im*(randn(20,20)-0.5))*sqrt(12).
\]

This matrix has real and complex entries that have mean $0$, variance
$1$ and kurtosis $\gamma=-1.2$.

An experimenter wants to understand how the smallest singular value compares
with that of the complex Gaussian matrix
\[
\tt randn(20,20)+im*randn(20,20).
\]

\begin{figure}

\includegraphics{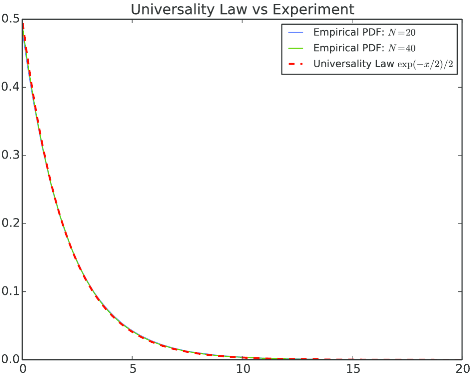}

\caption{Universality law vs. experiment: $n=20$ and $n=40$ already
resemble $n=\infty$.}\vspace*{-6pt}\label{fig1}
\end{figure}

The law for complex Gaussian matrices \cite
{Edelman,edelman1988eigenvalues} in this case valid for all finite
sized matrices, is that
$n\lambda_{\min}(AA^*)=n\sigma_{\min}^2(A)$ is exactly exponentially
distributed:
$f(x)=\frac{1}{2}e^{-x/2}$. Universality theorems say that the uniform
curve will match the Gaussian in the limit as matrix sizes go to
$\infty$. The experimenter obtains the curves in Figure~\ref{fig1} (taking both $n=20$
and $n=40$).

Impressed that $n=20$ and $n=40$ are so close, he or she might look at
the proof of the universality theorem only to find that no useful
bounds are available at $n=20,40$.

The results in this paper give the following correction in terms of the
kurtosis (when $\nu=0$):\vspace*{-3pt}
\[
f(x)= e^{-x/2} \biggl( \frac{1}{2} + \frac{\gamma}{n}\biggl(
\frac
{1}{4}-\frac
{x}{8}\biggr) \biggr) + O \biggl(\frac{1}{n^2}
\biggr).
\]

\begin{figure}

\includegraphics{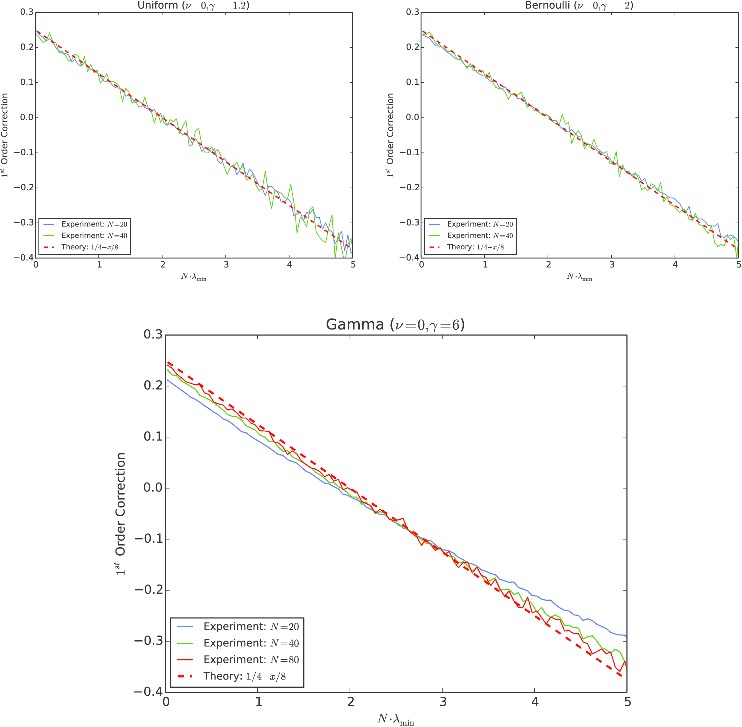}

\caption{Correction for square matrices Uniform, Bernoulli ($\nu=0$).
Monte Carlo simulations are histogrammed, 0th order term subtracted and
result multiplied by $n e^{x/2}/\gamma$. Bottom curve shows convergence
for $n=20,40,80$ for a distribution with positive kurtosis.}\label{fig2}
\end{figure}


On the bottom of Figure~\ref{fig1}, with the benefit of 60 computational
processors, we can magnify
the departure from universality with Monte Carlo experiments, showing
that the departure
truly fits $ \frac{\gamma}{n} (\frac{1}{4}-\frac{x}{8})e^{-x/2}$. This
experiment can be run and rerun
many times, with many distributions, kurtoses that are positive and
negative, small values of $n$,
and the correction term works very well.
Figure~\ref{fig2} shows
that the corrections converge as predicted for uniform, Bernoulli and Gamma distributed entries.\vspace*{-4pt}

\subsection{Example: $n+1$ by $n$ complex matrices (\texorpdfstring{$\nu=1$}{$nu=1$})}\label{sec24}

The correction to the density can be written as\vspace*{-2pt}
\[
f(x) = e^{-x/2} \biggl(\frac{1}{2} I_2 (s) +
\frac{1+\gamma
}{8n} \bigl(sI_1(s) -x I_2(s) \bigr)
\biggr)+ O\biggl(\frac{1}{n^2}\biggr),
\]
where $I_1(x)$ and $I_2(x)$ are Bessel functions and $s=\sqrt{2x}$.

\begin{figure}

\includegraphics{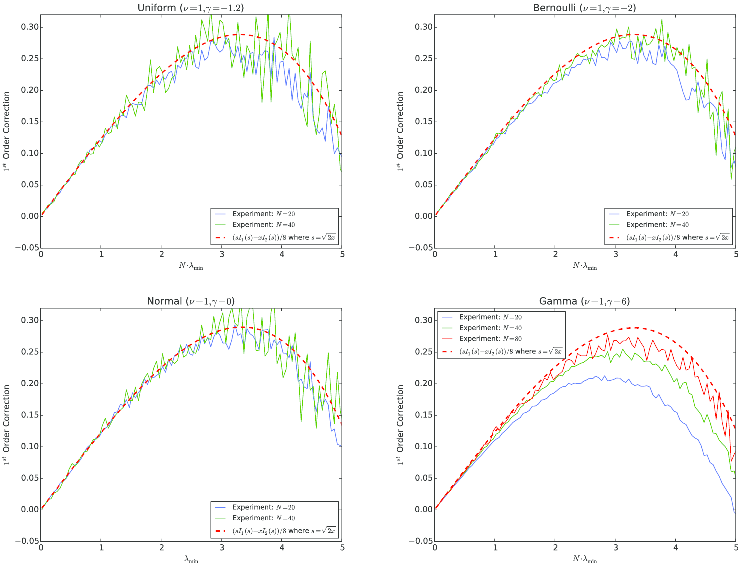}

\caption{Correction for $\nu=1$. Uniform, Bernoulli, Normal and Gamma;
Monte Carlo simulations are histogrammed, 0th order term subtracted and
result multiplied by $n e^{x/2}/(1+\gamma)$. Bottom right curve shows
convergence for $n=20,40,80$ for a distribution with positive kurtosis.}\label{fig3}
\end{figure}

Simulations are shown in Figure~\ref{fig3}.


\subsection{Example: $n+2$ by $n$ complex matrices (\texorpdfstring{$\nu=2$}{$nu=2$})}\label{sec25}
The correction to the density for $\nu=2$ can be written\vspace*{-2pt}
\begin{eqnarray*}
f(x) &=& \frac{1}{2}e^{-x/2} \biggl(\bigl[I^2_2(s)
- I_1(s)I_3(s)\bigr] \\[-3pt]
&&{}+\frac
{2+\gamma}{2n} \bigl[(x+4)
I_1^2 (s )-2 s I_0 (s ) I_1 (s
)-(x-2) I_2^2 (s ) \bigr] \biggr),
\end{eqnarray*}
where $I_0,I_1,I_2$, and $I_3$ are Bessel functions, and $s=\sqrt{2x}$.

\begin{figure}

\includegraphics{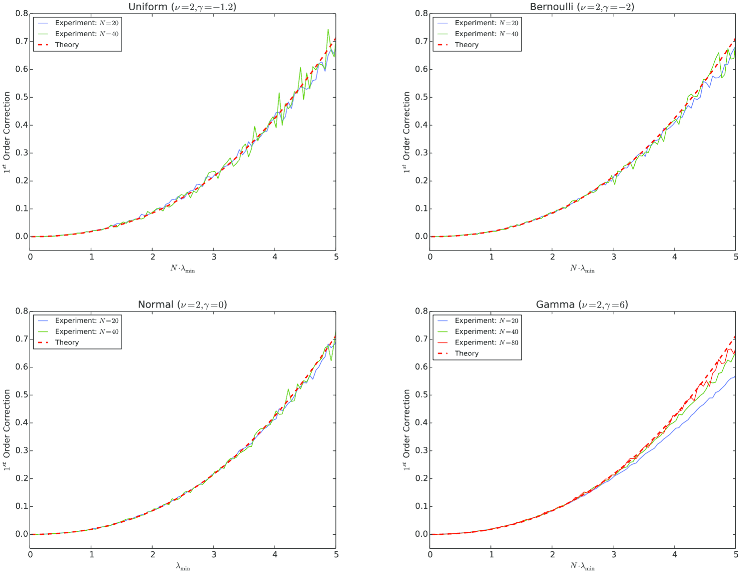}

\caption{Correction for $\nu=2$. Uniform, Bernoulli, Normal and Gamma;
Monte Carlo simulations are histogrammed, 0th order term subtracted and
result multiplied by $n e^{x/2}/(2+\gamma)$. Bottom right curve shows
convergence for $n=20,40,80$ for a distribution with positive kurtosis.}\label{fig4}
\end{figure}

Simulations are given in Figure~\ref{fig4}.

\section{Models and results}

In this section, we define the models we will study and state the results.
Let some real parameter $a>0$ be given.
Consider a matrix $M$ of size $p \times n$:
\[
M= W + aV,
\]
where:
\begin{itemize}
\item $V=(V_{ij})_{1\leq i\leq p; 1\leq j \leq n}$ has i.i.d. entries
with complex $\mathcal{N}_{\mathbb{C}}(0,1)$ distribution,
which means that both $\Re V_{ij}$ and $\Im V_{ij}$ are real i.i.d.
$\mathcal{N}(0, 1/2)$ random variables,


\item $W=(W_{ij})_{1\leq i\leq p; 1\leq j \leq n}$ is a random matrix
with entries being mutually independent random variables with
distribution $P_{ij}, 1\le j\le n$ independent of $n$ and $p$,
with uniformly bounded fourth moment,
\item $W$ is independent of $V$,
\item $\nu:=p-n\geq0$ is a fixed integer independent of $n$.
\end{itemize}

We then form the Gaussian divisible ensemble (also known as the
Johansson--Laguerre matrix):
%
\begin{equation}\label{eq1}
\frac{1}{n}M^*M= \biggl(\frac{1}{\sqrt n}(W+aV) \biggr)^* \biggl(
\frac
{1}{\sqrt n}(W+aV) \biggr).
\end{equation}
When $W$ is fixed, the above ensemble is known as the \emph{deformed
Laguerre ensemble.}

We assume that the probability distributions $P_{j,k}$ satisfy
%
%
\begin{equation}\label{4centrevariance}
\int z\,dP_{j,k}(z)=0,\qquad \int\bigl|zz^*\bigr| \,dP_{j,k}(z)=
\sigma^2_{\mathcal
C}=\frac{1}{4}.
\end{equation}
Here, the complex $\sigma^2_{\mathcal{C}} =2 \sigma^2_\Re$ represents
the complex variance.
Hypothesis (\ref{4centrevariance}) ensures the convergence of the
spectral measure of $\frac{1}{n}W^*W$ to the Marchenko--Pastur
distribution with density
%
\begin{equation}\label{densiteMPlimiteg1}
\rho_{\mathrm{PM}}(x)=\frac{2}{\pi}\frac{\sqrt{1-x}}{\sqrt x},\qquad
0\le x\le 1.
\end{equation}
Condition (\ref{4centrevariance}) implies also that the limiting
spectral measure of $\frac{1}{n}M^*M$ is then given by
Marchenko--Pastur's law with parameter $1/4+a^2$; we denote $\rho_a$
the density of this probability measure, that is, $\rho
_a(x)=(1+4a^2)^{-1/2}\rho_{\mathrm{PM}}(x/\break \sqrt{1+4a^2})$.

For technical reasons, we assume that the entries of $W$ have
subexponential tails: There exist $C,c,\theta>0$ so that for all
$j,\in
\mathbb{N}$, all $t\ge0$
\begin{equation}
\label{expdecay}
P_{j,k}\bigl(|z|\ge t\bigr)\le Ce^{-ct^\theta}.
\end{equation}
This hypothesis could be weakened to requiring enough finite moments.

Finally, we assume that the fourth moments do not depend on $j,k$ and
let $\kappa_4$ be the difference between the fourth moment of $P_{j,k}$
and the complex Gaussian case,
\[
\kappa_4=\int \bigl|zz^*\bigr|^2 \,dP_{j,k}-8^{-1}.
\]
(Thus, with the notation of Definition~\ref{defkurt}, $\kappa
_4=2\gamma
\sigma_{\Re}^4=2\kappa_4^{\Re}$.)

Then our main result is the following.
Let $\sigma:=\sqrt{4^{-1}+a^2}$ and for an Hermitian matrix $A$ denote
$\lambda_{\min}(A)=\lambda_1(A)\le\lambda_2(A)\le\cdots\le\lambda_n(A)$
the eigenvalues of $A$.

\begin{theorem}\label{maintheo2}
Let $F_n^\nu$ be the cumulative density function
of the hard edge of a Gaussian $p\times n$ matrix $V$, $\nu=p-n$, with
entries with complex variance $\sigma^2$:
\[
F_n^\nu(s)=\mathbb{P} \biggl(\sigma^2
\lambda_{\min}\bigl(VV^*\bigr)\leq\frac{
s}{n} \biggr).
\]
Then, for all $s> 0$, if our distribution has complex fourth cumulant
$\kappa_4=2\kappa_4^{\Re}$,
\[
\mathbb{P} \biggl( \lambda_{\min}\bigl(MM^*\bigr)\leq\frac{s}{n}
\biggr) =F_n^\nu (s) +\frac{s(F_n^\nu)'(s)}{ \sigma^4 n}
\kappa_4 +o\biggl(\frac{1}{n}\biggr).
\]
We note that this formula is scale invariant. It is equivalent to
\[
\mathbb{P} \biggl( \lambda_{\min}\bigl(MM^*\bigr)\geq\frac{s}{n}
\biggr) =1-F_n^\nu (s) +\frac{s(1-F_n^\nu)'(s)}{ \sigma^4 n}
\kappa_4 +o\biggl(\frac{1}{n}\biggr).
\]
\end{theorem}

Let $F_\infty^\nu(s)=\lim_{n\rightarrow\infty} F_n^\nu(s)$ be the limiting
cumulative distributive function in the Gaussian case.
$F^\nu_\infty$ is well known: we provide a few of its definitions in
Section~\ref{Gaussexp}; see also \eqref{forFinfty}. The difference
$F_n^\nu-F_\infty^\nu$
was derived by Schehr
\cite{Schehr} and Bornemann \cite{bornemann2015}, after it was
conjectured in the first version of this paper.
We can then deduce the following.
%
\begin{corollary}\label{tip}
For all integer number $\nu$,
\[
\mathbb{P} \biggl( \lambda_{\min}\bigl(MM^*\bigr)\leq\frac{s}{n}
\biggr) =F_\infty^\nu (s) + \biggl(\nu+\frac{\kappa_4}{\sigma^4}
\biggr) \frac{s(F_\infty^\nu)'(s)}{ n} +o\biggl(\frac
{1}{n}\biggr).
\]
\end{corollary}

\subsection*{Note}
Corollary~\ref{tip} is the convenient formulation we used for $\nu
=0,1,2$ in our showcase examples
in Sections~\ref{sec23}, \ref{sec24} and \ref{sec25}, respectively.

\subsection*{Note}
Corollary~\ref{tip} is remarkable because it states that the correction
term for $n$ being
a finite Gaussian as opposed to being infinite, and the correction term
for $n$ being
non-Gaussian as opposed to Gaussian ``line up,'' in that either way the
corrections
are multiples of $s(F^\nu_\infty)'(s)$. This could not be predicted by
Theorem~\ref{maintheo2} alone.

\subsection*{Note}
It is worth taking more of a close look between the formulation in
Theorem~\ref{maintheo2} and
Corollary~\ref{tip}. The first term in Theorem~\ref{maintheo2} is $n$
dependent, while in Corollary~\ref{tip}
the first term has reached its $n \rightarrow\infty$ limit.
Also of note is that the $F_n$ formulation in Theorem~\ref{maintheo2}
involves Laguerre polynomials (and exponentials.) The $F_\infty$
formulation in Corollary~\ref{tip} involves Bessel
functions (and exponentials).

For the Wigner ensemble, we consider the matrix
\[
M_n=\frac{1}{\sqrt{n}} (W+aV),
\]
where $W$ is a Wigner matrix with complex (resp., real) independent
entries above (resp., along) the diagonal $W_{ij}, 1\leq i\leq j \leq N
$ with law $P_{ij}$. We assume that the distributions $P_{ij}$
have subexponential moments: there exists $C,c>0$, and $\alpha>0$ such
that for all $t\ge0$ and all $1\leq i\leq j \leq N$
%
\begin{equation}
\label{expdecayWigner}
P_{ij}\bigl(|x|\ge t\bigr)\le C\exp\bigl\{ -ct^\alpha\bigr\},
\end{equation}
and satisfy
%
\begin{equation}
\label{momWigner}
\int x\,dP_{ij}(x)=0,\qquad \int|x|^2
\,dP_{ij}(x)=1/4,\qquad \int x^3\,dP_{ij}(x)=0.
\end{equation}
Again we assume that the fourth moments do not depend on $i,j$ and let
$\kappa_4$ be the difference between the fourth moment of $P_{ij}$
($j\neq i$) and the Gaussian case,
\[
\kappa_4=\int\bigl|zz^*\bigr|^2\,dP_{ij} -1/8.
\]
%
The other matrix $V$ is a GUE random matrix with i.i.d. $\mathcal
{N}_{\mathbb{C}
}(0,1)$ entries.
We denote by $\lambda_1 \leq\lambda_2 \leq\cdots\leq\lambda_n$ the
ordered eigenvalues of $M_n$.
By Wigner's theorem, it is known that the spectral measure of $M_n$
\[
\mu_n=\frac{1}{n}\sum_{i=1}^n
\delta_{\lambda_i}
\]
converges weakly to the semicircle distribution with density
%
\begin{equation}
\label{defsemicercle}
\sigma^{2\sigma}_{\mathrm{sc}}(x)= \frac
{1}{2\pi\sigma
^2}\sqrt{4
\sigma^2-x^2} \mathbh{1}_{|x|\leq2\sigma}; \qquad\sigma
^2=1/4+a^2.
\end{equation}

This is the Gaussian-divisible ensemble studied by Johansson \cite{Joh}.
We study the dependency of the one point correlation function $\rho_n$
of this ensemble, given as the probability measure on $\mathbb R$ so
that for any bounded measurable function $f$
\[
\mathbb E\Biggl[\frac{1}{n}\sum_{i=1}^n
f(\lambda_i)\Biggr]=\int f(x)\rho_n(x) \,dx
\]
as well as the
localization of the quantiles of $\rho_n$ with respect to the
quantiles of the limiting semicircle distribution. In particular, we
study the $1/n$ expansion of this localization, showing that it depends
on the fourth moment of $\mu$.
Define $N_n(x):= \frac{1}{n}\sharp\{i, \lambda_i \leq x\}$, and
$N_{\mathrm{sc}}(x)=\int_{-\infty}^x d\sigma_{\mathrm{sc}}(u)$, with $\sigma_{\mathrm{sc}}$
defined in (\ref{defsemicerclee}).
Let us define the quantiles $\hat\gamma_i$ (resp., $\gamma_i$) by
\[
\hat{\gamma}_i:=\inf \biggl\{y, \mathbb{E}N_n(y)=
\frac{i}{n} \biggr\}\qquad \mbox{respectively } N_{\mathrm{sc}}\bigl((-\infty,
\gamma_i]\bigr)=\frac{i}{n}.
\]
We shall prove the following.
%
\begin{theorem}\label{theotaovu} Let $\varepsilon>0$.
There exists a function $ D$ on $[-2+\varepsilon,2-\varepsilon]$,
independent of the distributions $P_{i,j}$, such that for all $x\in
[-2+\varepsilon,2-\varepsilon]$
\[
\rho_n(x)=\sigma_{\mathrm{sc}}(x)+\frac{1}{n}
\kappa_4 D(x)+o \biggl(\frac
{1}{n} \biggr).
\]
For all $i\in[n\varepsilon, n(1-\varepsilon)]$ for some $\varepsilon
>0$, one has that
%
\begin{eqnarray}
&&\hat\gamma_i-\gamma_i =\frac{\kappa_4}{2n} \bigl(2
\gamma_i^3-\gamma_i\bigr)+o\biggl(
\frac{1}{n}\biggr).
\end{eqnarray}
\end{theorem}
This is a version of the rescaled Tao--Vu conjecture 1.7 in \cite
{TaoVu} where $\mathbb E[\lambda_i]$ is replaced by $\hat\gamma_i$.
A similar result could be derived for Johansson--Laguerre ensembles. We
do not present the details of the computation here, which would
resemble the Wigner case. The function $D$ is computed explicitly in
Proposition~\ref{expprop}.

\section{Smallest singular values of \texorpdfstring{$n+\nu$}{$n+nu$} by $n$ complex Gaussian
matrices}\label{Gaussexp}
Theorem~\ref{maintheo2} depends on the partition function for Gaussian
matrices, which itself depends on $\nu$ and $n$.
In this section, we investigate these dependencies.

\subsection{Known exact results}

It is worthwhile to review what exact representations are known for the
smallest singular values of complex Gaussians.

We consider the finite $n$ density $f_n^\nu(x)$, the finite $n$
cumulative distribution $F_n^\nu(x)$ (we stress the $\nu$-dependency in
this section),
and their asymptotic values $f_\infty^\nu(x)$ and $F_\infty^\nu(x)$.
We have found the first form in the list below useful for symbolic and
numerical computation.
In the formulas to follow, we assume $\sigma^2_{\Re}=1$ so that a
command such as \verb+randn()+
can be used without modification for the real and imaginary parts.
All formulas concern $n\lambda_{\min}(AA^*)=n\sigma^2_{\min}(A)$ and
its asymptotics. We present in the array below eight different formulations
of the exact distribution $F_n^\nu$.

Some of these formulations allow one or both of $\nu$ or $n$ to extend
beyond integers to real positive values.
Assuming $\nu$ and $n$ are integers (\cite{Edelman}, Theorem~5.4),
the probability density $f_n^\nu(x)$
takes the form $x^\nu e^{-x/2}$ times a polynomial of degree $(n-1)\nu$
and $1-F_n^\nu(x)$ is $e^{-x/2}$ times a polynomial of degree $n\nu$.

\begin{rem*}
A helpful trick to compare normalizations used by different
authors is to inspect the
exponential term. The $2$ in $e^{-x/2}$ denotes total complex variance
$2$ (twice the real variance of 1).
In general the total complex variance $\sigma^2=2\sigma_{\Re}^2$ will
appear in the denominator.
\end{rem*}


In the next paragraphs, we discuss the eight formulations introduced above.

\subsubsection{Determinant: \texorpdfstring{$\nu$}{$nu$} by \texorpdfstring{$\nu$}{$nu$} determinant}

The quantities of primary use are the beautiful $\nu$ by $\nu$
determinant formulas for the distributions by Forrester and Hughes
\cite{Forrester1994}
in terms of Bessel functions and Laguerre polynomials.
The infinite formulas also appear in \cite{Forrester2010}, equation (8.98). Hereafter, $I_j$ denotes the modified Bessel
functions and $L_j$ the Laguerre polynomials:
\begin{eqnarray*}
F_\infty^\nu(x) &=& e^{-x/2} \det\bigl[I_{i-j}(\sqrt{2x}) \bigr]_{i,j=1,\ldots
,\nu}, \\
f_\infty^\nu(x) &=& \frac{1}{2}e^{-x/2} \det
\bigl[I_{2+i-j}(\sqrt{2x}) \bigr]_{i,j=1,\ldots,\nu}, \\
F_n^\nu(x) &=&  e^{-x/2}\det
\bigl[L_{n+i-j}^{(j-i)}(-x/2n) \bigr]_{i,j=1,\ldots,\nu}, \\
f_n^\nu(x) &=&
{
\biggl( \frac{x}{2n}  \biggr)}^{\nu}
\frac{(n-1)!}{ 2(n+\nu-1)!}
e^{-x/2}
\\
&&{}\times \det \bigl[L_{n-1+i-j}^{(j-i+2)} (-x/2n)
\bigr]_{i,j=1,\ldots,\nu}.
\end{eqnarray*}

Recall that $I_j(x)=I_{-j}(x)$. To facilitate reading of the relevant
$\nu$ by $\nu$ determinants, we provide expanded views:
\begin{eqnarray*}
&& \det\bigl[I_{i-j}(\sqrt{2x})\bigr]_{i,j=1,\ldots,\nu}\\
 &&\qquad = \left\vert
\matrix{
I_0 & I_1 &
I_2 & \cdots& I_{\nu-1}
\cr
I_1 & I_0 & I_1 & \cdots&
I_{\nu-2}
\cr
I_2 & I_1 & I_0 & \cdots& I_{\nu-3}
\cr
\vdots& \vdots& \vdots& \ddots& \vdots
\cr
I_{\nu-1} & I_{\nu-2} & I_{\nu-3} & \cdots&
I_0
}
\right\vert _{\mathrm{Bessel} \ \mathrm{functions} \ \mathrm{evaluated} \ \mathrm{at} \ \sqrt{2x}},
\\
&& \det\bigl[I_{2+i-j}(\sqrt{2x}) \bigr]_{i,j=1,\ldots,\nu}\\
&&\qquad=  \left\vert
\matrix{
I_2 & I_1 &
I_0 & \cdots& I_{\nu-3}
\cr
I_3 & I_2 & I_1 & \cdots& I_{\nu-4}
\cr
I_4 & I_1 & I_2 & \cdots& I_{\nu-5}
\cr
\vdots& \vdots& \vdots& \ddots& \vdots
\cr
I_{\nu+1} & I_{\nu} & I_{\nu-1} & \cdots&
I_2
}
\right\vert_{\mathrm{Bessel} \ \mathrm{functions} \ \mathrm{evaluated} \ \mathrm{at} \ \sqrt{2x}},
\\
&&\det \biggl[L_{n+i-j}^{(j-i)}\biggl(-\frac{x}{2n}\biggr)
\biggr]_{i,j=1,\ldots
,\nu} \\
&&\qquad= \left\vert
\matrix{
L_n & L_{n-1}^{(1)} & L_{n-2}^{(2)}
& \cdots& L_{n-\nu+1}^{(\nu-1)}
\hspace*{3pt}\cr
L_{n+1}^{(-1)} & L_n &
L_{n-1}^{(1)} & \cdots& L_{n-\nu+2}^{(\nu-2)}
\vspace*{3pt}\cr
L_{n+2}^{(-2)} & L_{n+1}^{(-1)}&
L_n &\cdots& L_{n-\nu+3}^{(\nu-3)}
\vspace*{3pt}\cr
\vdots& \vdots& \vdots& \ddots&\vdots
\vspace*{3pt}\cr
L_{n+\nu-1}^{(1-\nu)} & L_{n+\nu-2}^{(2-\nu)} &
L_{n+\nu
-3}^{(3-\nu
)} & \cdots& L_n
}
\right\vert_{\mathrm{evaluated} \ \mathrm{at} \ -x/{2n}},
\\
&&\det \biggl[L_{n-1+i-j}^{(j-i+2)} \biggl(-\frac{x}{2n}\biggr)
\biggr]_{i,j=1,\ldots,\nu} \\[3pt]
&&\qquad= \left\vert
\matrix{
L_{n-1}^{(2)} & L_{n-2}^{(3)} &
L_{n-3}^{(4)} & \cdots& L_{n-\nu}^{(\nu+1)}
\vspace*{3pt}\cr
L_n^{(1)} & L_{n-1}^{(2)}
& L_{n-2}^{(3)} & \cdots& L_{n-\nu+1}^{(\nu)}
\vspace*{3pt}\cr
L_{n+1} & L_n^{(1)} &
L_{n-1}^{(2)} & \cdots& L_{n-\nu+2}^{(\nu-1)}
\vspace*{3pt}\cr
\vdots& \vdots& \vdots& \ddots& \vdots
\vspace*{3pt}\cr
L_{n+\nu-2}^{(3-\nu)} & L_{n+\nu-3}^{(4-\nu)} &
L_{n+\nu-4}^{(5-\nu)} & \cdots& L_{n-1}^{(2)}
}
\right\vert_{\mathrm{evaluated} \  \mathrm{at} \ -x/{2n}}.
\end{eqnarray*}

The following Mathematica code symbolically computes these distributions:\vspace*{8pt}
\begin{verbatim}
   M[x_,v_]:= Table[BesselI[Abs[i-j],x],{i,v},{j,v}];
   m[x_,v_]:= Table[BesselI[Abs[2+i-j],x],{i,v},
              {j,v}];
 M[x_,n,v_]:= Table[LaguerreL[n+i-j,j-i,-x/(2*n)],
              {i,v},{j,v}];
m[x_,n_,v_]:= Table[LaguerreL[n-1+i-j,j-i+2,-x/(2*n)],
              {i,v},{j,v}];
  F[x_,v_ ]:= 1-Exp[-x/2]*Det [M[Sqrt[2 x],v]];
   f[x_,v_]:= (1/2)*Exp[-x/2]*Det[m[Sqrt[2 x],v]];
F[x_,n_,v_]:= 1-Exp[-x/2]*Det[M[x,n,v]];
f[x_,n_,v_]:= (x/(2n))^v*((n-1)!/(2(n+v-1)!))
              *Exp[-x/2]*Det[m[x,n,v]].
\end{verbatim}\vspace*{3pt}

\subsubsection{Painlev\'e III}


According to
\cite{Forrester2010}, equation (8.93),
\cite{bornemann2010}, pages 814--815,
\cite{Tw,tracywidom94},
we have the formula valid for all $\nu>0$:
\[
F_\infty^\nu(x)=1- \exp \biggl( -\int_0^{2x}
\sigma(s) \,\frac
{ds}{s} \biggr),
\]
where $\sigma(s)$ is the solution to a Painlev\'e III differential
equation. Please consult the references taking care to match the normalization.

\subsubsection{$n$ by $n$ determinant}

Following standard techniques to set up the multivariate integral and
applying a continuous version of the Cauchy--Binet theorem (Gram's formula)
\cite{mehta2004}, for example, Appendix~A.12  or \cite{tracyWidom98}, for example, equations (1.3)
and (5.2)  one can work out an
$n \times n$ determinant
valid for any $\nu$, so long as $n$ is an integer \cite{mike}:
\begin{eqnarray*}
F_n^{\nu}(x) &=& 1- \frac{\det(M(m,\nu,x/2)) } {\det(M(m,\nu,0))},\\[-8pt]
\
\end{eqnarray*}
where, if $\Gamma$ denotes the incomplete gamma function,
{\fontsize{10}{12}\selectfont
\begin{eqnarray*}
&&\! M(m,\nu,x)\\
&&\!\qquad = \pmatrix{
\Gamma(\nu+1,x) &
\Gamma(\nu+2,x) & \Gamma(\nu+3,x) & \cdots& \Gamma (\nu+m,x)
\cr
\Gamma(\nu+2,x) & \Gamma(\nu+3,x) & \Gamma(\nu+4,x) & \cdots& \Gamma (\nu+m+1,x)
\cr
\Gamma(\nu+3,x) & \Gamma(\nu+4,x) & \Gamma(\nu+5,x) & \cdots& \Gamma (\nu+m+2,x)
\cr
\vdots& \vdots& \vdots& \ddots& \vdots
\cr
\Gamma(\nu+m,x) & \Gamma(\nu+m+1,x) & \Gamma(\nu+m+2,x) & \cdots& \Gamma(
\nu+2m-1,x)}.
\end{eqnarray*}
}

\subsubsection{Remaining formulas in Table~\texorpdfstring{\protect\ref{tab1}}{3}}

The Fredholm determinant is a standard procedure. The multivariate
integral recurrence was computed in the real case in \cite{Edelman}
and in the complex case in \cite{Forrester1994}.
Various hypergeometric representations may be found in \cite
{dumitriu2003eigenvalue}, but to date we are not aware of the complex
representation of the confluent representation in \cite
{richards1982evaluation} which probably is worth pursuing.

\subsection{Asymptotics of smallest singular value densities of complex
Gaussians}\label{seconj}
A~very useful expansion extends a result from \cite{Forrester1994},
(3.29).

\begin{lemma}
\label{lmlag}
As $n\rightarrow\infty$, we have the first two terms in the asymptotic
expansion of scaled Laguerre polynomials
whose degree and constant parameter sum to~$n$:
\[
L_{n-k}^{(k)}(-x/n) \sim n^k \biggl\{
\frac{I_k(2\sqrt
{x})}{x^{k/2}}-\frac{1}{2n} \biggl( \frac{I_{k-2}(2\sqrt{x})}{x^{(k-2)/2}} \biggr) + O
\biggl(\frac{1}{n^2} \biggr) \biggr\}.
\]
\end{lemma}

\begin{pf}
We omit the tedious details but this (and indeed generalizations of
this result) may be computed either
through direct expansion of the Laguerre polynomial or through the
differential equation it satisfies.
\end{pf}

\begin{table}[t]
\caption{Exact results for smallest singular values of complex
Gaussians (smallest eigenvalues of complex Wishart or Laguerre ensembles)}
\label{tab1}
\begin{tabular*}{\textwidth}{@{\extracolsep{\fill}}lc@{}}
\hline
1. Determinant: $\nu$ by $\nu$ & \cite{Forrester1994,Forrester2010}\\
2. Painlev\'e III & \cite{Forrester2010}, equation (8.93) \\
3. Determinant: $n$ by $n$ & \cite{mike}\\
4. Fredholm determinant & \cite{tracywidom94,bornemann2010}\\
5. Multivariate integral recurrence & \cite{Edelman,Forrester1994} \\
6. Finite sum of schur polynomials (evaluated at $\mathbb{I}$) & \cite
{dumitriu2003eigenvalue}\\
7. Hypergeometric function of matrix argument & \cite{dumitriu2003eigenvalue}\\
8. Confluent hypergeometric function of matrix argument & \cite
{richards1982evaluation} \\
\hline
\end{tabular*}
\end{table}

We can use the lemma above to conjecture asymptotics of the
distribution $ F_n^{\nu}(x)$. This conjecture was recently proved independently
by Anthony Perret and Gr\'{e}gory Schehr
\cite{Schehr} and by Folkmar Bornemann \cite{bornemann2015}.
The first uses properties of the Jacobi matrix associated with modified Laguerre
polynomials and its implication to Painlev\'{e} while the second takes
a close look at the Fredholm determinant
and the asymptotics of Laguerre polynomials.
This result states as follows.

\begin{theorem}\label{conj1} Let $\nu$ be an integer number.
Let $F_n^{\nu}(x)$ be the distribution of $n\sigma_{\min}^2(A)$ of an
$n+\nu$ by $n$ complex Gaussian $A$.
We have the expansion
\[
F_n^{\nu}(x)=F_{\infty}^{\nu}(x) +
\frac{\nu}{2n}xf_{\infty}^{\nu
}(x) + O\biggl(\frac{1}{n^2}
\biggr).
\]
\end{theorem}

\subsection*{Note}
The above is readily checked to be scale invariant, so it is not
necessary to state the particular variances in the
matrix as long as they are equal.



\section{The hard edge of complex Gaussian divisible ensembles}
The hard edge denotes the location of the smallest eigenvalues of
sample covariance matrices when $\nu=p-n$ is a fixed integer.
\subsection{Reminder on Johansson--Laguerre ensemble}We here recall
some important facts about the Johansson--Laguerre ensemble, that we
use in the following.
\subsubsection*{Notation}
We\vspace*{1pt} call $\mu_{n,p}$ the law of the sample covariance matrix $\frac
{1}{n}M^*M$ defined in (\ref{eq1}). We denote by $\lambda_1 \leq
\lambda_2 \leq\cdots\leq\lambda_n$ the ordered\vspace*{1pt} eigenvalues of the
random sample covariance matrix $\frac{1}{n}M^*M$.
We also set
\[
H=\frac{W}{\sqrt n},
\]
and denote the distribution of the random matrix $H$ by $P_n$. The
ordered eigenvalues of
$HH^*$ are denoted by $y_1(H)\leq y_2(H)\leq\cdots\leq y_n(H)$.

We can now state the known results about the joint eigenvalue density
(j.e.d.) induced by the Johansson--Laguerre ensemble.
Propositions \ref{symprop}, \ref{symprop2} and \ref{CorrelHardEdge} are
derived in \cite{BAP}, Sections~3 and 7.
By construction, this is obtained as the integral w.r.t. $P_n$ of the
j.e.d. of the deformed Laguerre ensemble. We recall that the deformed
Laguerre ensemble denotes the distribution of the covariance matrix
$n^{-1} M M^*$ when $H$ is given. The latter has been first computed by
\cite{GW96} and \cite{JSV96}.

We now set
\[
s=\frac{a^2}{n}.
\]

\begin{proposition}\label{symprop}
The symmetrized eigenvalue measure on $\mathbb{R}_+^n$ induced by $\mu
_{n,p}$ has a density w.r.t. Lebesgue measure given by
%
\begin{equation}
\label{definitiondedQ} g (x_1,\ldots,x_n )=\int
dP_n(H) g \bigl(x_1,\ldots ,x_n;y(H) \bigr)
\end{equation}
with
\begin{eqnarray*}
&& g \bigl(x_1,\ldots,x_n;y(H) \bigr) \\
&&\qquad=\frac{\Delta(x)}{\Delta(y(H))}
\operatorname{det} \biggl( \frac{e^{-({y_i(H)+x_j})/({2t})}}{2t} I_{\nu} \biggl(\frac{\sqrt{y_i
(H)x_j}}{t}
\biggr) \biggl(\frac{x_j}{y_i(H)} \biggr)^{{\nu}/{2}} \biggr)_{i,j=1}^n,
\end{eqnarray*}
where $t=\frac{a^2}{2n}=\frac{s}{2}$, and $\Delta(x)=\prod_{i<j}(x_i-x_j)$.
\end{proposition}

From the above computation, all eigenvalue statistics can in principle
be computed.
In particular, the $m$-point correlation functions of $\mu_{n,p}$
defined by
\[
R_m(u_1, \ldots, u_m)=\frac{n!}{(n-m)!}
\int_{\mathbb
{R}_+^{n-m}}g(u_1, \ldots, u_n) \prod
_{i=m+1}^n \,du_i
\]
are given by the integral w.r.t. to $dP_n(H)$ of those of the deformed
Laguerre ensemble.
Let
\[
R_m \bigl(u_1,\ldots,u_m;y(H) \bigr)=
\frac{n!}{(n-m)!}\int_{\mathbb{R}
_+^{n-m}}g \bigl(u_1, \ldots,
u_n;y(H) \bigr) \prod_{i=m+1}^n
du_i
\]
be the $m$-point correlation function of the deformed Laguerre ensemble
(defined by the fixed matrix $H$). Then we have the following.
%
\begin{proposition}\label{symprop2}
\[
R_m (u_1,\ldots, u_m )=\int
_{M_{p,n}(\mathbb
{C})}\,dP_n(H)R_m \bigl(u_1,
\ldots, u_m;y(H) \bigr).
\]
In particular,
\[
P\biggl(\lambda_{\min}\biggl(\frac{MM^*}{n}\biggr)\ge a\biggr)=\int
_a^\infty R_1(u) \,du= \int
_a^\infty\!\int_{M_{p,n}(\mathbb{C})}dP_n(H)R_1
\bigl(u;y(H) \bigr) \,du.
\]
\end{proposition}
The second remarkable fact is that the deformed Laguerre ensemble
induces a determinantal random point field, that is all the $m$-point
correlation functions are given by the determinant of a $m\times m$
matrix involving the same \emph{correlation kernel}.
%
\begin{theorem}\label{CorrelHardEdge}
Let $m$ be a given integer. Then one has that
\[
R_m \bigl(u_1,\ldots, u_m;y(H) \bigr)=\det
\bigl( K_n \bigl(u_i,u_j; y(H) \bigr)
\bigr)_{i,j=1}^m,
\]
where the correlation kernel $K_n$ is defined by
\begin{eqnarray*}
&& K_n \bigl(u,v; y(H) \bigr)\\
&&\qquad= \frac{e^{\nu i \pi}}{i \pi s^3} \int
_{\Gamma
}\!\int_{\gamma}\,dw \,dz wzK_B
\biggl(\frac{2z\sqrt u}{s},\frac{2w\sqrt
v}{s} \biggr) \biggl(\frac{w}{z}
\biggr)^{\nu} \exp\biggl\{\frac
{w^2-z^2}{s}\biggr\}
\\
&&\qquad\quad{}\times\prod_{i=1}^n \frac{w^2-y_i(H)}{z^2-y_i(H)}
\Biggl(1-s\sum_{i=1}^n\frac{y_i(H)}{(w^2-y_i(H))(z^2-y_i(H))}
\Biggr),
\end{eqnarray*}
where the contour $\Gamma$ is symmetric around 0 and encircles the
$\pm
\sqrt{y_i(H)}$, $\gamma$ is the imaginary axis oriented positively $0
\longrightarrow+\infty$, $0 \longrightarrow-\infty$, and $K_B$ is the
kernel defined by
%
\begin{equation}
\label{defBesK} K_B(x,y)=\frac{x I_{\nu
}^{\prime}(x)I_{\nu
}(y)-y I_{\nu}^{\prime}(y)I_{\nu}(x)}{x^2-y^2}.
\end{equation}
\end{theorem}
There are two important facts about this determinantal structure.
The fundamental characteristic of the correlation kernel is that it
depends only on the spectrum of $HH^*$ and more precisely on its
spectral measure.
Since we are interested in the determinant of matrices with entries
$K_n (x_i,x_j;y(H) )$, we can consider the correlation kernel
up to a conjugation: $K_n (x_i,x_j;\break y(H) ) \frac
{f(x_i)}{f(x_j)}$. This has no impact on correlation functions and we
may use this fact later.

For ease of exposition, we drop from now on the dependency of the
correlation kernel $K_n$ on the spectrum of $H$ and write $K_n(u,v)$ for
$K_n (u,v; y(H) )$.
The goal of this section is to deduce Theorem~\ref{maintheo2} by a
careful asymptotic analysis of the above formulas.
Set
%
\begin{equation}
\label{defialpha}
\alpha=\sigma^2/4,
\end{equation}
with $\sigma=\sqrt{1/4+a^2}$.

We recall that it was proved in \cite{BAP} that
%
\begin{equation}
\label{forFinfty}
\lim_{n \to\infty} \mathbb{P} \biggl(
\lambda_{\min} \biggl(\frac{MM^*}{n} \biggr)\geq \frac
{\alpha s}{n^2}
\biggr)=\det (I-\tilde{K}_B )_{L^2(0, s)},
\end{equation}
where $\tilde{K}_B$ is the usual Bessel kernel
%
\begin{equation}
\label{deftildeK}\tilde{K}_B(u,v):=e^{\nu i \pi
}K_B(i
\sqrt{u}, i\sqrt{v})
\end{equation}
with $K_B$ defined in \eqref{defBesK}. This is a universality result as
the limiting distribution function is the same as that of the smallest
eigenvalue of the Laguerre ensemble (see, e.g., \cite{Forrester}).

\subsection{Asymptotic expansion of the partition function at the hard
edge}
\label{subsecasan}
The main result of this section is to prove the following expansion for
the partition function at the hard edge: recall that $\alpha$ is given
by (\ref{defialpha}).
%
\begin{theorem}\label{teethe}
Assume that the distributions $P_{jk}$
satisfy the
assumptions (\ref{expdecay}) and (\ref{4centrevariance}). Then
there exists a nonnegative function $g_n^0$, depending on $n$, so that
\[
\mathbb{P} \biggl( \lambda_{\min} \biggl(\frac{MM^*}{n} \biggr)\geq
\frac
{\alpha s}{n^2} \biggr) =g_n^0(s) +\frac{1}{n}
\,\partial_\beta g_n^0(\beta s)\Big|_{\beta=1}
\int dP_n(H)\bigl[\Delta_n(H)\bigr]+o \biggl(
\frac
{1}{n} \biggr),
\]
where
\[
\Delta_n(H)=\frac{-1}{v_c^\pm m_{\mathrm{MP}}'(v_c^+)} X_n\bigl(v_c^+
\bigr)
\]
with $X_n(z) =\sum_{i=1}^n \frac{1}{y_i(H)-z}-nm_{\mathrm{MP}}(z)$,
$m_{\mathrm{MP}}(z)$ is the Stieltjes transform of the
Marchenko--Pastur distribution $\rho_{\mathrm{MP}}$, $(y_i(H))_{1\le i\le
n}$ are the eigenvalues of $H$,
and $v_c^+=(w_c^\pm)^2$ where
%
\begin{equation}
\label{crit}
w_c^\pm=\pm i (R-1/R)/2, \qquad R:=
\sqrt{1+4a^2}.
\end{equation}
\end{theorem}
We will estimate the term $\int dP_n(H)[\Delta_n(H)]$ in terms of the
kurtosis in the next section.
%
\begin{rem}
The function $g_n^0$ is universal, in the sense that it
does not
depend on the detail of the distributions $P_{jk}$.
\end{rem}
%
\subsubsection{Expansion of the correlation kernel}
Let $z_c^\pm$ be the critical points of
%
\begin{equation}
\label{defFn}E_n(w):=w^2/a^2+
\frac{1}{n}\sum_{i=1}^n \ln
\bigl(w^2-y_i(H)\bigr),
\end{equation}
where the $y_i(H)$ are the eigenvalues of $H^*H$. Then we have the
following lemma. Let $K_n$ be the kernel defined in Theorem~\ref
{CorrelHardEdge}.

\begin{lemma} \label{expo} There exists a smooth function $A$ which is
independent of $\kappa_4$ such that for all $u,v$
\begin{eqnarray*}
&&\frac{\alpha}{n^2}K_n \bigl(u \alpha n^{-2}, v\alpha
n^{-2};y(H) \bigr)\\
&&\qquad=
\tilde{K}_B(u, v)+
\frac{A(u,v)}{n}+\biggl(\biggl(\frac
{z^+_c}{w^+_c}\biggr)^2-1\biggr)
\,\frac{\partial}{\partial\beta} \Big|_{\beta
=1}\beta\tilde{K}_B(\beta u,
\beta v)+o \biggl(\frac{1}{n} \biggr),
\end{eqnarray*}
for $H$ in a set with probability greater than $1- e^{-n^{1/2}}$, and
where
$\tilde{K}_B$ has been defined in (\ref{deftildeK}). Note that
$\frac{z^+_c}{w^+_c}=\frac{z^-_c}{w^-_c}$.
\end{lemma}
\begin{pf}
To focus on local eigenvalue statistics at the hard edge, we consider
\[
u = \biggl(\frac{a^2}{2n r_0} \biggr) ^2 x; \qquad v= \biggl(
\frac{a^2}{2n
r_0} \biggr)^2y \qquad \mbox{where $r_0$ will be
fixed later.}
\]
As $\nu=p-n$ is a fixed integer independent of $n$, this readily
implies that the Bessel kernel shall not play a role in the large
exponential term of the correlation kernel.
In other words, the large exponential term to be considered is $E_n$
defined in (\ref{defFn}).
The correlation kernel can then be rewritten as
%
\begin{eqnarray}
K_n(u,v) &=& \frac{1}{i \pi s^3} e^{\nu i \pi} \int
_{\Gamma}\!\int_{\gamma
} dw \,dz w z
K_B \biggl(\frac{zx^{1/2}}{r_0},\frac{wy^{1/2}}{r_0} \biggr) \biggl(
\frac{w}{z} \biggr)^{\nu}
\nonumber
\\[-8pt]
\\[-8pt]
\nonumber
&&{}\times\exp{\bigl\{n
E_n(w)-nE_n(z)\bigr\}}\tilde{G}(w,z),
\end{eqnarray}
where
\begin{eqnarray*}
\tilde{G}(w,z)&:=&a^2G(w,z)=1-s\sum_{i=1}^n
\frac
{y_i(H)}{(w^2-y_i(H))(z^2-y_i(H))}
\\
&=& \frac{a^2}{2} \frac{wE_n'(w)-zE_n'(z)}{w^2-z^2}.
\end{eqnarray*}
We note that $E_n(w)=H_n(w^2)$ where $H_n(w)=w/a^2+\frac{1}{n}\sum_{i=1}^n \ln(w-y_i(H))$.

We may compare the exponential term $E_n$ to its ``limit,'' using the
convergence of the spectral measure of $H^*H=\frac{1}{n} W^*W$ to the
Marchenko--Pastur distribution $\rho_{\mathrm{MP}}$.
Set
\[
E(w):= w^2/a^2+\int\ln\bigl(w^2-y\bigr) \,d
\rho_{\mathrm{MP}}(y).
\]
It was proved in \cite{BAP} that this term has two conjugated critical points
satisfying $E'(w)=0$ and are given by $w_c^{\pm}$ defined in \eqref{crit}.
Let us also denote by $z_c^{\pm}$ the true nonreal critical points
(which can be seen to exist and be conjugate \cite{BAP}) associated to
$E_n$. These critical points do depend on $n$ but for ease of notation
we do not stress this dependence. These critical points satisfy
\[
E_n'\bigl(z_c^{\pm}\bigr)=0,
\qquad z_c^+=-z_c^-
\]
and it is not difficult to see that they are also on the imaginary axis.

We now refer to the results established in \cite{BAP} to claim the
following facts:
\begin{itemize}
\item there exists a constant $C$ so that for
any $\xi \in (0,1)$
\begin{equation}
\label{remarque17}
\bigl|z_c^{\pm}-w_c^{\pm}\bigr|
\leq Cn^{-\xi}
\end{equation}
with probability greater than $1-e^{-n^{2-2\xi}}$
for $n$ large enough. In the
sequel we will take $\xi = 3/4$. This comes from concentration results for
the spectral measure of $H$ established in \cite{GZ} and \cite{Bai2}  and formula (\ref{dep}).
%
\item Fix $\theta>0$. By the saddle point analysis performed in \cite
{BAP}, the contribution of the parts of the contours $\gamma$ and
$\Gamma$ within $ \{|w-z_c^{\pm}|\geq n^{\theta}n^{-1/2}\}$ is $O(e^{-c
n^{\theta}})$ for some $c >0$. This contribution ``far from the critical
points'' is thus exponentially negligible. In the sequel, we will
choose $\theta=1/11$. The choice of $1/11$ is arbitrary.
\item We can thus restrict both the $w$ and $z$ integrals to
neighborhoods of width $n^{1/11}n^{-1/2}$ of the critical points
$z_c^{\pm}$.
\end{itemize}
Also, we can assume that the parts of the contours $\Gamma$ and
$\gamma
$ that will contribute to the asymptotics are symmetric w.r.t.
$z_c^{\pm
}$. This comes from the fact that the initial contours exhibit this
symmetry and from the location of the critical points.
A plot of the oriented contours close to critical points is given in
Figure~\ref{fig5}.
%
\begin{figure}

\includegraphics{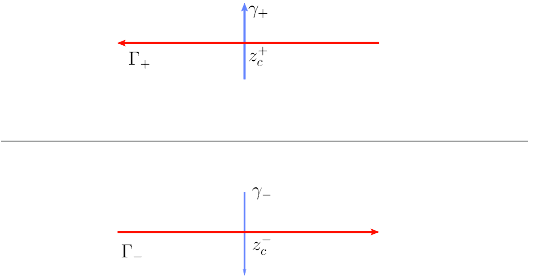}

\caption{Contours close to the critical points.}\label{fig5}\vspace*{-6pt}
\end{figure}

Let us now make the change of variables
\[
w= z_c^{1}+ sn^{-1/2}; \qquad z=
z_c^{2}+ tn^{-1/2},
\]
where $z_c^{1}, z_c^2$ are equal to $z_c^{+}$ or $z_c^-$ depending on
the part of the contours $\gamma$ and $\Gamma$ under consideration and
$s,t$ satisfy $|s|, |t| \leq n^{1/11}$.
Then we perform the Taylor expansion of each of the terms arising in
both $z$ and $w$ integrands. Then one has that
%
\begin{eqnarray}
&& e^{n E_n(z_c^{\pm}+ sn^{-1/2})-nE_n(z_c^{\pm})}
\nonumber\\[-2pt]
&& \qquad =e^{ E_n''(z_c^{\pm}){s^2}/{2}+ \sum_{i=3}^5E_n^{(i)}(z_c^{\pm
}){s^i}/({i!n^{i/2-1}})}\bigl(1+O\bigl(n^{-23/22}\bigr)
\bigr)
\nonumber\\[-9pt]
\\[-9pt]
&&\qquad =e^{ E_n''(z_c^{\pm}){s^2}/{2}}+ \frac{1}{n^{1/2}}\underbrace {e^{ E_n''(z_c^{\pm}){s^2}/{2}}
\frac{E_n^{(3)}(z_c^{\pm
})}{6}s}_{e_1(s)}
\nonumber
\\[-2pt]
\nonumber
&&\qquad\quad{}+ \frac{1}{n}
\underbrace{e^{ E_n''(z_c^{\pm}){s^2}/{2}} \biggl(\frac{E_n^{(4)}(z_c^{\pm}) s^4}{4!}+ \biggl(\frac{E_n^{(3)}(z_c^{\pm
})}{6}
\biggr)^2 \frac{s^6}{2} \biggr)}_{e_2(s)}+ o \biggl(
\frac
{1}{n} \biggr)e^{ E_n''(z_c^{\pm}){s^2}/{2}},
\end{eqnarray}
as $|s|\leq n^{1/11}$.
For each term in the integrand, one has to consider the contribution of
equal or opposite critical points. In the following, we denote by $z_c,
z_c^1, z_c^2$ any of the two critical points (allowing $z_c$ to take
different values with a slight abuse of notation).
We then perform the Taylor expansion of each of the functions arising
in the integrands. This yields the following four expansions:
%
\begin{equation}
wz=z_c^1z_c^2+
n^{-1/2}\underbrace{\bigl( sz^2_c+t
z^1_c\bigr)}_{v_1(s,t)}+\frac
{1}{n}
\underbrace{st}_{v_2(s,t)},
\end{equation}
and
%
\begin{eqnarray}
&& g \biggl(z_c^1+ \frac{s}{n^{1/2}},
z_c^2+ \frac
{t}{n^{1/2}} \biggr)\nonumber\\
&&\qquad=\frac{E_n''(z_c)}{2}
\mathbh{1}_{z_c^1=z_c^2} +\frac{1}{\sqrt{n}}\underbrace{ \biggl(s
\frac{\partial}{\partial x_1}+t \frac{\partial}{\partial x_2} \biggr)G(x_1,
x_2) \Big|_{z_c^1,
z_c^2}}_{g_1(s,t)}
\\
&&\qquad\quad {}+\frac{1}{n}
\underbrace{ \biggl( \biggl(\frac{s^2} 2\frac{\partial
^2}{\partial x_1^2}G(x_1,
x_2)+\frac{t^2} 2 \frac{\partial
^2}{\partial
x_2^2}+st \frac{\partial^2}{\partial x_2\,\partial x_1}
\biggr)G(x_1, x_2) \Big|_{z_c^1, z_c^2}
\biggr)}_{g_2(s,t)}\nonumber\\[-2pt]
&&\qquad\quad{}+o \biggl(\frac{1}{n} \biggr).\nonumber
\end{eqnarray}
One also has\vspace*{-2pt}
%
\begin{eqnarray}
\biggl(\frac{w}{z} \biggr)^{\nu
} &=& \bigl(z_c^1/z_c^2
\bigr)^{\nu}+ n^{-1/2}\underbrace{\bigl(z_c^1/z_c^2
\bigr)^{\nu}\biggl(\frac{\nu s}{z_c^1}-\frac
{\nu
t}{z_c^2}
\biggr)}_{r_1(s,t)}
\nonumber
\\[-9pt]
\\[-9pt]
\nonumber
&&\qquad{}+\frac{1}{n}\underbrace {
\bigl(z_c^1/z_c^2
\bigr)^{\nu
} \biggl( \frac{\nu(\nu-1)s^2}{(z_c^1)^2} +\frac{\nu(\nu+1)t^2}{(z_c^2)^2}-
\frac{\nu^2st}{z_c^1
z_c^2} \biggr)}_{r_2(s,t)}+\,o \biggl(\frac{1}{n} \biggr).
\end{eqnarray}
Last, one has that
%
\begin{eqnarray}
&& \!\!\! K_B \biggl(\frac{zx^{1/2}}{r_0},\frac
{wy^{1/2}}{r_0}
\biggr) \nonumber\\
&&\!\!\!\qquad= K_B \biggl(\frac{z_cx^{1/2}}{r_0},\frac{z_c
y^{1/2}}{r_0} \biggr)
\nonumber
\\[-8pt]
\\[-8pt]
\nonumber
&&\!\!\!\quad\qquad{}+ \frac{1}{\sqrt n}\underbrace{ \biggl( s\frac{\partial}{\partial x_1
}+t
\frac{\partial}{\partial x_2} \biggr) \Big|_{z_c,z_c} K_B \biggl(
\frac
{x_1 x^{1/2}}{r_0},\frac{x_2y^{1/2}}{r_0} \biggr)}_{h_1(s,t)}
\\
&&\!\!\!\quad\qquad{}+
\frac{1}{n}\underbrace{ \biggl( \frac{s^2}{2}\frac{\partial
^2}{\partial x_1^2 }+
\frac{t^2}{2}\frac{\partial^2}{\partial x_2^2 }+st \frac{\partial^2}{\partial x_1 \,\partial x_2 } \biggr)
\Big|_{z_c,z_c}\! K_B \biggl(\frac{x_1 x^{1/2}}{r_0},\frac{x_2y^{1/2}}{r_0}
\biggr)}_{h_2(s,t)}+\,o \biggl(\frac{1}{n} \biggr).\hspace*{-25pt}\nonumber
\end{eqnarray}
In all the lines above, $z_c^1/z_c^2=\pm1$ as critical points are
either equal or opposite. Also one can note that the $o$ are uniform as
long as $|s|, |t|<n^{1/11}$.

We now choose
\[
r_0=\bigl|w_c^+\bigr|.
\]
Combining the whole contribution of neighborhoods of a pair of equal
critical points, for example, denoted by $K_n(u,v)_{equal}$, we find
that it has an expansion of the form
%
\begin{eqnarray}
&&\frac{a^4}{4n^2 r_0^2}K_n(u,v)_{\mathrm{equal}}\nonumber\\
&&\qquad= \sum
_{z_c=z_c^{\pm}}
\frac{ \pm}{4i \pi}e^{\nu i \pi}\int
_{\mathbb{R}}\!\int_{i
\mathbb{R}}\,ds \,dt \frac
{|z_c|^2}{r_0^2}
\Biggl( K_B\biggl(\frac{z_{c}x^{1/2}}{|w_c^+|},\frac
{z_cy^{1/2}}{|w_c^+|}\biggr) \\
&&\qquad\quad{}+
\sum_{i=1}^2\frac{h_i(s,t)}{n^{i/2}}+o \biggl(
\frac{1}{n} \biggr) \Biggr)
\nonumber\\
&& \qquad\quad{}\times \Biggl( \frac
{E_n''(z_c)}{2}+\sum
_{i=1}^2\frac{g_i(s,t)}{n^{i/2}}+o \biggl(
\frac
{1}{n} \biggr) \Biggr) \Biggl(1+\sum_{i=1}^2
\frac
{r_i(s,t)}{n^{i/2}}+o \biggl(\frac{1}{n} \biggr) \Biggr)\nonumber
\\
&&\qquad\quad{}\times
\Biggl(1+\sum_{i=1}^2v_i(s,t)n^{-i/2}
z_c^{-2}\Biggr) \biggl(\exp{\bigl\{ E_n''(z_c)
\bigl(s^2-t^2\bigr)/2\bigr\}} \biggl(1+o \biggl(
\frac{1}{n} \biggr)
\nonumber\\
&& \qquad\quad{}+ n^{-1/2}\bigl(e_1(s)-e_1(t)
\bigr)+\frac{1}{n} \bigl(-e_1(s)e_1(t)+e_2(s)-e_2(t)
\bigr) \biggr) \biggr),\nonumber
\end{eqnarray}
where $h_i, e_i, r_i, v_i$ and $g_i$ defined above have no singularity.

It is not difficult also to see that $h_1, g_1, r_1, e_1$ are odd
functions in $s$ as well as in $t$: because of the symmetry of the
contour, their contribution will thus vanish. The first nonzero lower
order term in the asymptotic expansion will thus come from the combined
contributions $h_1g_1, g_1 r_1, r_1 h_1, h_1 e_1, g_1e_1,\break r_1e_1 ,
r_1v_1,\ldots$ and those from $h_2, g_2, r_2, e_2, v_2$.
Therefore, one can check that one gets the expansion
\begin{eqnarray*}
&&\frac{\alpha}{n^2}K_n \biggl(\frac{\alpha x}{ n^{2}},\frac{\alpha y}{
n^{2}}
\biggr)_{\mathrm{equal}}\\
&&\qquad=\frac{e^{i\nu\pi}}{2} \biggl( \frac
{|z_c^{\pm
}|}{|w_c^{\pm}|}
\biggr)^2 K_B \biggl(\frac{z_{c}^{\pm
}x^{1/2}}{|w_c^+|},
\frac{z_c^{\pm}y^{1/2}}{|w_c^+|} \biggr)+\frac{a_1
(z_c^{\pm}; x,y)}{n}+o \biggl(\frac{1}{n}
\biggr),
\end{eqnarray*}
where $a_1$ is a function of $z_c^{\pm}, x,y$ only. $a_1$ is a smooth
and nonvanishing function a priori.

We can write the first term above as $ (\frac{z_c^\pm}{w_c^\pm
} )^2 \tilde{K}_B ( (\frac{z_c^\pm}{w_c^\pm
}
)^2 x,  (\frac{z_c^\pm}{w_c^\pm} )^2y )$ so that we
deduce that
\begin{eqnarray*}
&&e^{i\nu\pi} \biggl( \frac{z_c^{\pm}}{w_c^{\pm
}} \biggr)^2
K_B \biggl(\frac{z_{c}^{\pm}x^{1/2}}{|w_c^+|},\frac
{z_c^{\pm}y^{1/2}}{|w_c^+|} \biggr)
\\
&&\qquad =
\tilde{K}_B(x,y)+ \biggl( \biggl( \frac{z_c^{\pm}}{w_c^{\pm
}}
\biggr)^2- 1 \biggr)\partial_\beta \bigl(\beta
\tilde{K}_B (\beta x,\beta y)\bigr)\Big|_{\beta=1}+o
\bigl(z_c^\pm -w_c^\pm \bigr).
\end{eqnarray*}
One can do the same thing for the combined contribution of opposite
critical points and get a similar result. We refer to \cite{BAP} for
more detail about this fact.
Summing these terms yield a contribution of order $\frac{1}{n}$.
However, it is clear that, as $a_1$ is smooth and using (\ref{remarque17}),
%
\begin{equation}
\label{A}
a_1\bigl(z_c^{\pm};x,y
\bigr)=a_1\bigl(w_c^{\pm};x,y\bigr)+o(1).
\end{equation}
Note that $w_c^{\pm}$ does not depend on the exact distributions
$P_{jk}$, but only on the limiting Marchenko--Pastur distribution $\rho
_{\mathrm{MP}}$. As a consequence,
there is no fourth moment contribution in this $\frac{1}{n}$ terms. We
denote the contribution of the deterministic error from all the
combined (equal or not) critical points by $A(x,y)/n$. This completes
the proof of the lemma.
\end{pf}

\subsubsection{Asymptotic expansion of the density}

The distribution of the smallest eigenvalue of $M_n$ is defined by
%
\[
\mathbb{P} \biggl( \lambda_{\min} \biggl(\frac{MM^*}{n} \biggr)\geq
\frac
{\alpha s}{n^2} \biggr) =\int dP_n(H)\det\bigl(I-\tilde{K}_n\bigl(y(H)\bigr)\bigr)_{L^2(0,s)},
\]
for $H$ in a set with overwhelming\vspace*{1pt} probability and where $\tilde{K}_n$ is the rescaled correlation kernel
$\frac{\alpha}{n^2}K_n (x \alpha n^{-2}, y\alpha n^{-2}; y(H))$. In the above we choose $\alpha=(a^2/2r_0)^2$. The limiting correlation kernel is
then, at the first order, the Bessel kernel:
\[
\tilde{K}_B(x,y):=e^{\nu i \pi}K_B(i\sqrt{x},
i\sqrt{y}).
\]

Hereafter we drop the dependency in $y(H)$ to simplify the notations.
The error terms are ordered according to their order of magnitude:
the first-order error term, in the order of $O(n^{-1})$, can thus come
from two terms, namely:
\begin{itemize}
\item[--] The deterministic part that is $A(x,y)/n$, which we have seen is
independent of~$\kappa_4$.
\item[--]  The kernel (arising 4 times due to the combination of critical points)
%
\begin{eqnarray}
&& e^{\nu i \pi}\biggl(\frac
{z_c^{+}}{|w_c^+|}\biggr)^2K_B
\biggl(\frac
{z_c^{+}}{|w_c^+|} (\sqrt{x},\sqrt{y}) \biggr)
\nonumber
\\[-8pt]
\label{B}
\\[-8pt]
\nonumber
&&\qquad= \tilde{K}_B(x,y)+
\int_1^{(z_c^+/w_c^+)^2}\frac{\partial}{\partial\beta}\beta \tilde{K}_B(\beta x, \beta y) \,d\beta.
\end{eqnarray}
\end{itemize}
Lemma~\ref{expo} and the arguments above (\ref{A}), (\ref{B}) give the
following:
\begin{eqnarray*}
&&\frac{\alpha}{n^2}K_n \bigl(x \alpha n^{-2}, y\alpha
n^{-2} \bigr)
\cr
&& \qquad = \tilde{K}_B(x, y)+
\frac
{A(x,y)}{n}+{\bigl(\bigl(z_c^+/w_c^+
\bigr)^2-1\bigr)}\frac
{\partial}{\partial\beta} \Big|_{\beta=1}\beta\tilde{K}_B(\beta x, \beta y)+o \biggl(\frac{1}{n} \biggr).
\end{eqnarray*}
We insist that the kernel $A$ is universal in the sense that it does
not depend on the detail of the distributions $P_{jk}$.

The Fredholm determinant can be developed to obtain that
%
\begin{eqnarray}
&&\det(I-\tilde{K}_n)_{L^2(0,s)}\nonumber
\\
&&\qquad=\sum_{k}\frac{(-1)^k}{k!}\int
_{[0,s]^k}\det \bigl( \tilde{K}_n (x_i,
x_j) \bigr)_{i,j=1}^k\prod
dx_i
\nonumber
\\[-8pt]
\\[-8pt]
\nonumber
&&\qquad=\sum_{k}
\frac{(-1)^k}{k!}\int_{[0,s]^k}\det \bigl( \tilde{K}_B
(x_i, x_j) \bigr)_{i,j=1}^k \det
\bigl(I +G(x_i,x_j) \bigr)_{i,j=1}^k
\prod dx_i\\
&&\qquad\quad {}+o\biggl(\frac{1}{n}\biggr),\nonumber
\end{eqnarray}
where we have set
\[
G(x_i,x_j)= \bigl( \tilde{K}_B
(x_i, x_j) \bigr)_{1\leq i,j\leq
k}^{-1}
\bigl(B(x_i, x_j) \bigr)_{i,j=1}^k
\]
with
%
\begin{equation}\label{f13}
\qquad B(x_i,x_j)=\frac{A(x_i,x_j)}{n}+ \biggl( \biggl(
\frac
{z_c^+}{w_c^+} \biggr)^2-1 \biggr)\frac{\partial}{\partial\beta
}
\Big|_{\beta=1}\beta\tilde{K}_B(\beta x_i, \beta
x_j)+o \biggl(\frac
{1}{n} \biggr).\hspace*{-20pt}
\end{equation}
The matrix $ ( \tilde{K}_B (x_i, x_j) )_{i,j=1}^k$ is
indeed invertible for any $k$.

Therefore, up to an error term in the order $o (\frac{1}{n}
)$ at most,
\begin{eqnarray*}
&&\det(I-\tilde{K}_n)_{L^2(0,s)}\\
&&\qquad=\det(I-\tilde{K}_B)
\\
&&\qquad\quad{}+\sum
_{k}\frac{(-1)^k}{k!}\int_{[0,s]^k}\det
\bigl( \tilde{K}_B (x_i, x_j)
\bigr)_{i,j=1}^k \operatorname{Tr} \bigl( G(x_i,x_j)
\bigr)_{i,j=1}^k\, dx +o\biggl(\frac{1}{n}\biggr).
\end{eqnarray*}
Now if we just consider the term which is linear in
$((z_c^+/w^+_c)^2-1)$ which will bring the contribution depending on
the fourth cumulant,
we have that the correction is
\begin{eqnarray*}
L &:=& \sum_{k}\frac{(-1)^k}{k!}\int
_{[0,s]^k}\det \bigl( \tilde{K}_B
(x_i, x_j) \bigr)_{i,j=1}^k \\[-2pt]
&&{}\times\operatorname{Tr}\bigl( \tilde{K}_B^{-1} \,\partial_\beta
\beta\tilde{K}_B (\beta x_i, \beta
x_j) \bigr)_{i,j=1}^k \,dx|_{\beta=1}
\\[-2pt]
&=& \partial_\beta\sum_{k}
\frac{(-1)^k}{k!}\int_{[0,s]^k}\det \bigl( \tilde{K}_B
(x_i, x_j) \bigr)_{i,j=1}^k \\[-2pt]
&&{}\times\operatorname{Tr}\bigl( \log\beta\tilde{K}_B (\beta x_i, \beta
x_j) \bigr)_{i,j=1}^k \,dx|_{\beta=1}.
\end{eqnarray*}
As $\tilde{K}_B$ is trace class, we can write
\[
\operatorname{Tr}\bigl( \log\beta\tilde{K}_B (\beta x_i,
\beta x_j)\bigr)_{i,j=1}^k =\log \det \bigl(
\beta\tilde{K}_B (\beta x_i,\beta
x_j) \bigr)_{i,j=1}^k.
\]
Therefore, we have
%
\begin{eqnarray}
L &=&\partial_\beta\sum_{k}
\frac{(-1)^k}{k!}\int_{[0,s]^k}\det \bigl( \tilde{K}_B
( x_i, x_j) \bigr)_{i,j=1}^k\nonumber \\[-2pt]
&&{}\times\log
\det \bigl( \beta \tilde{K}_B (\beta x_i,\beta
x_j) \bigr)_{i,j=1}^k \,dx|_{\beta
=1}
\nonumber
\\[-2pt]
\label{f14}
&=&\partial_\beta\sum_{k}
\frac{(-1)^k}{k!}\int_{[0,s]^k}\det \bigl( \beta
\tilde{K}_B (\beta x_i, \beta x_j)
\bigr)_{i,j=1}^k \,dx|_{\beta=1}
\\[-2pt]
\nonumber
&=&\partial_\beta\sum_{k}
\frac{(-1)^k}{k!}\int_{[0,s \beta
]^k}\det \bigl( \tilde{K}_B
( y_i, y_j) \bigr)_{i,j=1}^k
\,dy_i|_{\beta
=1}
\\[-2pt]
\nonumber
&=&\partial_\beta\det(I-\tilde{K}_B)_{L^2(0, s\beta)}
|_{\beta=1}.
\end{eqnarray}
Hence, since $ \det(I-\tilde{K}_B)_{L^2(0, s\beta)}$ is the
leading order in the expansion of $ \mathbb{P} ( \lambda
_{\min}(\frac
{MM^*}{n})\geq\frac{\alpha s}{n^2} ) $
plugging \eqref{f13} into \eqref{f14} shows that there exists a
function $g_n^0$ [whose leading order is $\det(I-\tilde{K}_B)_{L^2(0, s\beta)}$] so that
%
\begin{eqnarray}
&& \mathbb{P} \biggl( \lambda_{\min}\biggl(
\frac{MM^*}{n}\biggr)\geq\frac{\alpha
s}{n^2} \biggr)
\nonumber
\\[-8pt]
\label{exp1}
\\[-8pt]
\nonumber
&&\qquad =g_n^0(s)
+ \partial_\beta g_n^0(\beta
s)|_{\beta=1} \int dP_n(H) \biggl[ \biggl(\frac{z_c^+}{w_c^+}
\biggr)^2-1 \biggr]+o \biggl(\frac{1}{n} \biggr).
\end{eqnarray}
%

\subsubsection{An estimate for \texorpdfstring{$(\frac{z_c^+}{w_c^+} )^2-1$}{$({z_c^+}/{w_c^+})^2-1$}}
Let
\[
X_n(z)=\sum_{i=1}^n
\frac{1}{y_i(H)-z} -n m_{\mathrm{MP}}(z),
\]
where $z\in\mathbb C\setminus\mathbb{R}_+$.
Let us express $(z_c^+)^2-(w_c^+)^2$ in terms of $X_n$.
The\vspace*{1pt} critical point $z_c^+$ of $E_n$ lies in a neighborhood of the
critical point $w_c^+$
of $E$. So $u_c^+=(z_c^+)^2$ is in a neighborhood of
$v_c^+=(w_c^+)^2$. These points are the solutions with nonnegative
imaginary part of
\[
\frac{1}{a^2} +\frac{1}{n}\sum_{i=1}^n
\frac{1}{u_c^+-y_i(H)}=0, \qquad \frac{1}{a^2} +\int\frac{1}{v_c^+-y}\,d
\rho_{\mathrm{MP}}(y)=0.
\]
Therefore, it is easy to check that
\[
-\int\frac{u_c^+-v_c^+}{(v_c^+-y)^2}\,d\rho_{\mathrm{MP}}(y)+\frac{1}{n}
X_n\bigl(v_c^+\bigr)= o \biggl(\frac{1}{n},
\bigl(z_c^+-w_c^+\bigr) \biggr)
\]
which gives
%
\begin{equation}
\label{dep} \biggl(\frac{z_c^+}{w_c^+} \biggr)^2-1=
\frac
{1}{v_c^+ m_{\mathrm{MP}}'(v_c^+)}\frac{1}{n} X_n\bigl(v_c^+
\bigr) +o \biggl(\frac
{1}{n} \biggr).
\end{equation}
The proof of Theorem~\ref{teethe} is therefore complete. In the next
section, we estimate the expectation of $X_n(v_c^+)$ to get the
correction in \eqref{exp1}.
\subsection{The role of the fourth moment}
In this section, we compute\break $\mathbb E[X_n(v_c^+) ]$, which with
Theorem~\ref{teethe}, will allow to complete the proof of Theorem~\ref{maintheo2}.
\subsubsection{The expected value $\mathbb{E}[X_n(v_c^+) ]$}
\label{subsecCLT}

In this section, we give the asymptotics of the mean of $X_n(z)$. Such
type of estimates is now well known,
and can for instance be found in Bai and Silverstein book \cite{BaiSil}
for either Wigner matrices or Wishart matrices with $\kappa_4=0$. We
refer to \cite{BaiSil}, Theorem~9.10,  for a precise statement.
In the more complicated setting of $F$-matrices, we refer the reader to
\cite{Shu}. In the case where $\kappa_4\neq0$, the asymptotics of the
mean have been computed in \cite{Najim} and \cite{bai2010}.

\begin{proposition}\label{asympmean}
Let $z\in\mathbb{C}\setminus
\mathbb{R}_+ \cap\{Z \in\mathbb{C}, \Im Z
\geq0\}$. Then
\[
\lim_{n\rightarrow\infty} {\mathbb E}\bigl[X_n(z)\bigr]=A(z)-
\kappa_4 B(z)
\]
with $A$ independent of $\kappa_4$, and if $m_{\mathrm{MP}}(z)=\int
(x-z)^{-1} \,d\rho_{\mathrm{MP}}(x)$,
%
\begin{equation}
\label{fB}
B(z)=\frac{m_{\mathrm{MP}}(z)^2}{(1+
({m_{\mathrm{MP}}(z))/4})^2(z+({zm_{\mathrm{MP}}(z))/2})}.
\end{equation}
\end{proposition}

Note that\vspace*{1.5pt} the above result follows from a simple expansion (up to the
$1/N$ order) of the normalized trace of the resolvent $\frac
{1}{n}\operatorname{Tr} (\frac{WW^*}{n} -zI )^{-1}$ for a complex number~$z$
with nonzero imaginary part. We recall \cite{PM}
that
\[
m_{\mathrm{MP}}(z):=\lim_{n \to\infty}\frac{1}{n}\operatorname{Tr}
\biggl(\frac
{WW^*}{n} -zI \biggr)^{-1}
\]
is uniquely defined as the solution with nonnegative imaginary part
of the equation
%
\begin{equation}
\label{eqm}
\frac{1}{1+{1}/{4}m_{\mathrm{MP}}(z)}=-zm_{\mathrm{MP}}(z).
\end{equation}

\subsubsection{Estimate at the critical point}
Since $m_{\mathrm{MP}}(v_c^+)=a^{-2}$ and $v_c^+<0$ so that Proposition~\ref
{asympmean} applies, we deduce from \eqref{fB} that
there exists a constant $c(v_c^+)$ independent of $\kappa_4$ such that
\begin{eqnarray*}
\label{111}
&&\mathbb{E}\bigl[X_n\bigl(v_c^+\bigr)
\bigr]=c\bigl(v_c^+\bigr) -\kappa_4 \frac{a^{-4}}{(1+({1}/{4})
a^{-2})^2}
\frac{1}{v_c^+(1+({1}/{2})a^{-2})} +o(1).
\end{eqnarray*}
Moreover, we have
%
\begin{equation}\label{1111}
v_c^+=-\frac{a^4}{{1}/{4}+a^2}=-\frac
{4a^4}{1+4a^2}
\end{equation}
and by \eqref{eqm}, after taking the derivative, we find
\[
m_{\mathrm{MP}}'(z)=-{ \frac{m_{\mathrm{MP}}(z) (1+m_{\mathrm{MP}}(z)/4)}{z(1+m_{\mathrm{MP}}(z)/2)}},
\]
so that
at the critical point we get
%
\begin{eqnarray}
m_{\mathrm{MP}}'\bigl(v_c^+\bigr) &=&
\frac{(4 a^2+1)^2}{16
a^6(a^2+{1}/{2})},
\nonumber
\\[-8pt]
\label{11111}
\\[-8pt]
\nonumber
{v_c^+} m_{\mathrm{MP}}'
\bigl(v_c^+\bigr) &=& - \frac{(1+4a^2)}{4 ({1}/{2}+a^2)} = -\frac{a^{-2}(1+{1}/({4a^2}))}{1+{1}/({2a^2})}.
\end{eqnarray}

Therefore, with the notation of Theorem~\ref{teethe}, and using (\ref
{111}), (\ref{1111}) and \eqref{11111}, we find constants $C$
independent of $\kappa_4$ (and which may change from line to line)
so that
\begin{eqnarray*}
&&\int dP_n(H)\bigl[\Delta_n(H)\bigr]\\
&&\qquad=
\frac{1}{v_c^+ m_{\mathrm{MP}}'(v_c^+)} \mathbb E\bigl[n\bigl(m_n\bigl(v_c^+
\bigr)-m_{\mathrm{MP}}\bigl(v_c^+\bigr)\bigr)\bigr]+o(1)
\\
&&\qquad=
\frac{1+{1}/({2a^2})}{a^{-2}(1+{1}/({4a^2}))} \kappa_4\frac{a^{-4}}{(1+{1}/({4a^2}))^2}\frac{1}{1+{1}/({2a^2})}
\frac
{1+4a^2}{4a^4}+C+o(1)
\\
&&\qquad= \frac{16 \kappa_4}{(1+4a^2)^2}+C+o(1)=\frac{\kappa_4}{ ({1}/{4}+a^2)^2}+C+o(1).
\end{eqnarray*}

Rescale the matrix $M$ by dividing it by $\sigma$ so as to standardize
the entries. Combining Theorem~\ref{teethe} and the above, we have
therefore found that the deviation of the smallest eigenvalue are such that
\[
\mathbb{P} \biggl(\lambda_{\mathrm{min}} \biggl(\frac{MM^*}{n\sigma^2} \biggr)\ge
\frac{s}{n^2} \biggr)=\mathbf{g}_n(s)+\frac{\gamma}{2n} s
\mathbf {g}_n'(s) +o \biggl(\frac{1}{n} \biggr),
\]
where $\gamma$ is the kurtosis defined in Definition \eqref{defkurt}.
At this point $\mathbf{g}_n$ is identified
to be the distribution function at the hard edge of the Laguerre
ensemble with variance~$1$, as it corresponds to the case where $\gamma
=0$. Theorem~\ref{maintheo2} follows.

\section{The bulk of Gaussian divisible ensembles}
We here choose to consider the deformed GUE instead of the deformed
Laguerre ensemble. Indeed, while the arguments are completely similar,
the technicalities in the deformed Laguerre ensemble are more involved.
To ease the reading, we here present the simplest ensemble.

\subsection{Deformed GUE in the bulk}

Let $W=(W_{ij})_{i,j=1}^n$ be a Hermitian\vspace*{1pt} Wigner matrix of size $n$.
The entries $W_{ij}$ $1\leq i<j\leq n$ are i.i.d. with distribution
$P_{ij}$. The entries along the diagonal are i.i.d. real random
variables with law $P_{ii}$ independent of the off diagonal entries.
We assume that $P_{ij}, P_{ii}, 1\leq i\leq j \leq N$ have
subexponential tails and satisfy (\ref{expdecayWigner}) and (\ref
{momWigner}). The fourth moment of the $P_{ij}$'s is also assumed not
to depend on $i,j$.
Let also $V$ be a GUE random matrix with i.i.d. $\mathcal{N}_{\mathbb{C}
}(0,1)$ entries and consider the rescaled matrix
\[
M_n=\frac{1}{\sqrt n} (W+aV).
\]
We denote by $\lambda_1 \leq\lambda_2 \leq\cdots\leq\lambda_n$ the
ordered eigenvalues of $M_n$.
By Wigner's theorem, it is well known that the spectral measure of $M_n$
\[
\mu_n=\frac{1}{n}\sum_{i=1}^n
\delta_{\lambda_i}
\]
converges weakly to the semicircle distribution with density
%
\begin{equation}
\label{defsemicerclee}\sigma_{\mathrm{sc}}^{2\sigma}(x)= \frac
{1}{2\pi\sigma
^2}\sqrt{4
\sigma^2-x^2} 1_{|x|\leq2\sigma}; \qquad
\sigma^2=1/4+a^2.
\end{equation}

This is the deformed GUE ensemble studied by Johansson \cite{Joh}.
In this section, we study the localization of the eigenvalues $\lambda
_i$ with respect to the quantiles of the limiting semicircle
distribution. We study the $\frac{1}{n}$ expansion of this
localization, showing that it depends on the fourth moment of $P_{ij}$,
and prove Theorem~\ref{theotaovu}.

The route we follow is similar to that we took in the
previous section for Wishart matrices: we first obtain a $\frac{1}{n}$
expansion of the correlation functions of the Deformed GUE. The
dependency of this expansion in the fourth moment of $P_{ij}$ is then derived.

\subsection{Asymptotic analysis of the correlation functions}
Let $\rho_n$ be the one point correlation function of the Deformed GUE.
We prove in this subsection the following result, with $z_c^\pm
,w_c^\pm$
critical points similar to those of the last section, which we will
define precisely in the proof.
%
\begin{proposition}\label{expprop}
For all $\varepsilon>0$, uniformly on $u\in[-2\sigma+\varepsilon
,2\sigma-\varepsilon]$, we have
\[
\rho_n(u)=\sigma_{\mathrm{sc}}^{2\sigma}(u) + \mathbb E
\biggl[ \biggl(\frac{\Im
z_c^+(u)}{\Im w_c^+(u)}-1 \biggr)\biggr] \sigma_{\mathrm{sc}}^{2\sigma}(u)
+o \biggl(\frac{1}{n} \biggr),
\]
where $z_c^+$ depends
on the eigenvalues of $W$.
\label{proprhon}
\end{proposition}
%
\begin{pf}
Denote by $y_1(\frac{W}{\sqrt{n}})\leq y_2 (\frac{W}{\sqrt{n}})\leq
\cdots\leq y_n(\frac{W}{\sqrt{n}})$ the ordered eigenvalues of
$W/\sqrt n$.
Johansson
\cite{Joh}, (2.20) and (2.21) (see also \cite{BrezHik}) proves that,
for a fixed $W/\sqrt n$, the eigenvalue density of $M_n$ induces a
determinantal process with correlation kernel given by
\begin{eqnarray*}
&& K_n \biggl(u,v; y\biggl(\frac{W}{\sqrt{n}}\biggr) \biggr)\\
&&\qquad=
\frac{n}{(2i \pi)^2 } \int_{\Gamma} \,dz \int_{\gamma}\,dw
e^{n (E_v(w)-E_v(z))} \frac{1-e^{({(u-v)zn})/{a^2}}}{z(u-v)} g_n(z,w),
\end{eqnarray*}
where
\[
E_v(z)=\frac{(z-v)^2}{2a^2}+\frac{1}{n}\sum
_{i=1}^n \ln \biggl(z-y_i\biggl(
\frac{W}{\sqrt{n}}\biggr) \biggr),
\]
and
\[
g_n(z,w)=E_u'(z)+z\frac{E_u'(z)-E_u'(w)}{z-w}.
\]
The contour $\Gamma$ has to encircle all the $y_i$'s and $\gamma$ is
parallel to the imaginary axis.

We now consider the asymptotics of the correlation kernel in the bulk,
that is close to some point $u_0 \in(-2\sigma+\delta, 2\sigma
-\delta)$
for some $\delta>0$ (small). We recall that we can consider the
correlation kernel up to conjugation: this follows from the fact that
$\det ( K_n (x_i,x_j;y(\frac{W}{\sqrt{n}}) )
)=\det ( K_n (x_i,x_j;y(\frac{W}{\sqrt{n}}) )\frac
{h(x_i)}{h(x_j)} )$, for any nonvanishing function $h$.
We omit some details in the next asymptotic analysis as it closely
follows the arguments of \cite{Joh} and those of Section~\ref{subsecasan}.

Let then $u,v$ be points in the bulk with
\[
u=u_0+\frac{\alpha x}{n},\hspace*{18pt} v =u_0+\frac{\alpha\tilde{x}}{n};\hspace*{18pt}
u_0=\sqrt {1+4a^2} \cos(\theta_0),\hspace*{18pt}
\theta_0 \in(2\varepsilon, \pi-2\varepsilon).
\]
The constant $\alpha$ will be fixed afterward.
Then the approximate large exponential term to lead the asymptotic
analysis is given by
\[
\tilde{E}_{v} (z)=\frac{(z-v)^2}{2a^2}+\int\ln(z-y) \,d\sigma
_{\mathrm{sc}}^{1}(y) .
\]
In the following, we note $R_0=\sqrt{1+4a^2}=2\sigma$.

We recall the following facts from \cite{Joh}, Section~3. Let
$u_0=\sqrt
{1+4a^2} \cos(\theta_0)$ be a given point in the bulk.
\begin{itemize}
\item The approximate critical points, that is the solutions of $\tilde{E}_{u_0} '(z)=0$ are given by
\[
w_c^{\pm}(u_0)= \biggl(R_0e^{i \theta_c}
\pm\frac{1}{R_0e^{i \theta
_c}} \biggr)\Big/2.
\]
The true critical points satisfy $E_{u_0}'(z)=0$. Among the solutions,
we disregard the $n-1$ real solutions which are interlaced with the
eigenvalues $y_1(\frac{W}{\sqrt{n}}), \ldots, y_n(\frac{W}{\sqrt{n}})$.
The two remaining solutions are complex conjugate with nonzero
imaginary part and we denote them by
$z_c^{\pm}(u_0)$.
Furthermore \cite{GZ} and~\cite{Bai} prove that
\[
\bigl|z_c^{\pm}-w_c^{\pm}\bigr| \leq
Cn^{-\xi}
\]
with probability greater than $1-e^{-n^{2-2\xi}}$
for $n$ large enough and any point $u_0$ in the bulk of the
spectrum. In the sequel we will take $\xi = 3/4$.
\item We now fix the contours for the saddle point analysis. The steep
descent/ascent contours can be chosen as
\begin{eqnarray*}
\gamma &=&z_c^+(v) +it,\qquad t\in\mathbb{R},
\\
\Gamma &=& \bigl\lbrace
z_c^{\pm}(r), r=R_0 \cos(\theta), \theta\in (
\varepsilon, \pi-\varepsilon) \bigr\rbrace\cup \bigl\lbrace z_c^{\pm
}
\bigl(R_0\cos(\varepsilon)\bigr)+x, x>0 \bigr\rbrace
\\
&&{}\cup \bigl
\lbrace z_c^{\pm}\bigl(-R_0\cos(\varepsilon)
\bigr)-x, x>0 \bigr\rbrace.
\end{eqnarray*}
It is an easy computation [using that $\Re E_{u_0}'' (w)>0$ along
$\gamma$] to check that the contribution of the contour $\gamma\cap\{
|w-z_c^{\pm}(v)| \geq n^{1/12-1/2}\}$ is exponentially negligible.
Indeed there exists a constant $c>0$ such that
\[
\biggl| \int_{\gamma\cap\{ |w-z_c^{\pm}(v)| \geq n^{1/12-1/2}\}
}e^{n\Re
(E_{u_0}(w)-E_{u_0}(z_c^+(v)))} \,dw \biggr| \leq e^{- cn^{1/6}}.
\]
Similarly, the contribution of the contour $\Gamma\cap\{|w-z_c^{\pm
}(v)| \geq n^{1/12-1/2}\}$ is of order $e^{- cn^{1/6}}$ that of a
neighborhood of $z_c^{\pm}(v)$.
\end{itemize}

For ease of notation, we now denote $z_c(v):=z_c^+(v)$.
We now modify slightly the contours so as to make the contours
symmetric around $z_c^{\pm}(v)$. To this\vspace*{1pt} aim, we modify the $\Gamma$
contour as follows: in a neighborhood of width $n^{1/12-1/2}$
we replace $\Gamma$ by a straight line through $z_c^{\pm}(v)$ with
slope $z_c'(v)$. This slope is well defined as
\[
z_c'(v)=\frac{1}{E_v'' (z_c(v))}\neq0,
\]
using that $|z_c^{\pm}(v)-w_c^{\pm}(u_0)|\leq n^{-\xi}$.
We refer to Figure~\ref{figdeformation}, to define the new contour
$\Gamma'$ which is more explanatory.
\begin{figure}

\includegraphics{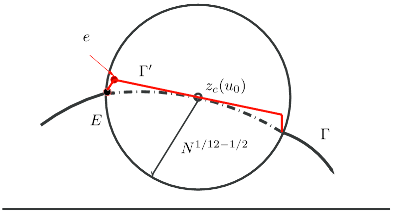}

\caption{Modification of the $\Gamma$ contour.}
\label{figdeformation}
\end{figure}
Denote by $E$ the leftmost point of $\Gamma\cap\{w,
|w-z_c(v)|=n^{1/12-1/2}\}$. Then there exists $v_1$ such that
$E=z_c(v_1)$. We then define $e$ by $e=z_c(v)+z_c'(v) (v_1-v)$.
We then draw the segment $[e, z_c(v)]$ and draw also its symmetric to
the right of $z_c(v)$.
Then it is an easy fact that
\[
|E-e|\leq C n^{2(1/12-1/2)}\qquad \mbox{for some constant } C.
\]
Here, we have used that $e, E$ both lie within a distance
$n^{1/12-1/2}$ from $z_c(v)$. It follows that
\[
\forall z \in[e, E], \qquad \bigl| \Re \bigl( n E_v(z)-n
E_v(E) \bigr) \bigr| \leq Cn n^{3 (1/12-1/2)} \ll  n^{1/6}.
\]
This follows from the fact that $|E_v'(z)|=O(n^{1/12-1/2})$ along the
segment $[e, E]$.
This is now enough as $\Re n E_v(E)> \Re nE_v(z_c)+cn^{1/6}$ to ensure
that the deformation has no impact on the asymptotic analysis.

We now make the change of variables $z=z_c^{\pm}(v)+\frac{t}{\sqrt n},
w=z_c^{\pm}(v)+\frac{s}{\sqrt n}$
where $|s|, |t| \leq n^{1/12}$.
We examine the contributions of the different terms in the integrand.
We first consider $g_n$. We start with the combined contribution of
equal critical points, for example, $z$ and $w$ close to the same
critical point. In this case we have, noting that
$E'_u(z)-E'_v(z)=a^{-2}(v-u)$
that
\begin{eqnarray*}
z^{-1}g_n(w,z)&=& E''_v
\bigl(z_c(v)\bigr) + \frac{\alpha(\tilde{x}-x)}{a^2nz_c(v)}+\frac{1}{\sqrt{n}}\biggl(
\frac{E_v^{(3)}(z_c(v))}{2 }(s+t)+ \frac{E_u'' (z_c(v))}{z_c(v)
}t\biggr)
\\
&&{}+\frac{1}{n}
\biggl( \frac{E_v^{(4)}(z_c(v))}{3! }\bigl(s^2+t^2+st\bigr)+
\frac
{E_v^{(3)}(z_c(v))}{2 z_c(v)} t^2- \frac{E_v'' (z_c(v))}{z_c(v)^2 } t^2 \biggr)
\\
&&{}+o \biggl(\frac{1}{n} \biggr).
\end{eqnarray*}
%
On the other hand, when $w$ and $z$ lie in the neighborhood of
different critical points, one gets that
\begin{eqnarray*}
z^{-1} g_n(w,z) &=&\frac{\alpha(\tilde{x}-x)}{a^2nz_c^{\mp}(v)}+
E_u'' \bigl(z_c^\mp(v)
\bigr)\frac{t}{z_c^\mp(v)\sqrt n }
\cr
&&{}+\frac{E_v^{(2)}(z_c^{\pm}(v))t-E_v^{(2)}(z_c^{\mp}(v))s}{
(z_c^{\mp
}-z_c^{\pm}) \sqrt n}+O \biggl(
\frac{1}{n} \biggr),
\end{eqnarray*}
%
where the $O (\frac{1}{n} )$ depends on the second and third
derivative of $E_v$ only.

We next turn to the second term, which depends on $z$ only. One has that
\[
1-e^{(x-\tilde{x})\alpha a^{-2}z_c^{\pm}}=1- e^{(x-\tilde{x})
\alpha
a^{-2}\Re z_c^+}e^{i \pm(x-\tilde{x})\alpha a^{-2}\Im z_c}.
\]
We then perform the same Taylor expansion as in Section~\ref{subsecasan} of all the terms in the integrands. As the contours are symmetric
around $z_c(u_0)$, the first nonzero term in the expansion is in the
scale of $\frac{1}{n}$.
Furthermore, apart from constants, one has that
\begin{eqnarray*}
\frac{\alpha}{n}K_n \biggl(u,v;y\biggl(\frac{W}{\sqrt{n}}\biggr)
\biggr)&=&\frac
{e^{(x-\tilde{x}) ({\alpha}/{a^2})\Re z_c^+}}{2 i\pi(x-\tilde{x})} \bigl(e^{i (x-\tilde{x})({\alpha}/{a^2})\Im z_c^+}-e^{-i
(x-\tilde{x})({\alpha}/{a^2})\Im z_c^+} \bigr)
\\
&& {}+\frac{C(x,\tilde{x})}{n} +o \biggl(\frac{1}{n} \biggr).
\end{eqnarray*}
The function $C(x,\tilde{x})$ does not depend on the detail of the
distributions of the entries of $W$.
We now choose $\alpha=\sigma_{\mathrm{sc}}^{2\sigma}(u_0)^{-1}$ where $\sigma
_{\mathrm{sc}}^{2\sigma}$ is the density of the semicircle distribution
defined in (\ref{defsemicerclee}). It has been proved in \cite{Joh} that
$\Im w_c^+(u_0)=\pi a^2\sigma_{\mathrm{sc}}^{2\sigma}(u_0)$.
Setting then
\[
\beta:=\Im z_c^+(u_0)/\Im\bigl(w_c^+(u_0)
\bigr)
\]
we then obtain that
\[
\frac{\alpha}{n}K_n \biggl(u,v;y\biggl(\frac{W}{\sqrt{n}}\biggr)
\biggr) e^{-(x-\tilde{x}) ({\alpha}/{a^2})\Re z_c^+}=\frac{\sin\pi
\beta
(x-\tilde{x})}{\pi(x-\tilde{x})}+\frac{C'(x,\tilde{x})}{n}+o \biggl(
\frac{1}{n} \biggr).
\]
The constant $C'(x,\tilde{x})$ does not depend on the distribution of
the entries of $W$. This proves Proposition~\ref{proprhon} since
\begin{eqnarray*}
\rho_n(x)&=&\mathbb{E}\biggl[\frac{1}{n} K_n
\biggl(u,u;y\biggl(\frac{W}{\sqrt
{n}}\biggr) \biggr)\biggr]
\\
&=&\frac{1}{\alpha} \mathbb E[\beta] + \frac{C'(x,x)}{n}+ o \biggl(
\frac
{1}{n} \biggr)
\\
&=&\sigma_{\mathrm{sc}}^{2\sigma}(u_0)+ \sigma_{\mathrm{sc}}^{2\sigma}(u_0)
\mathbb E\bigl[(\beta-1)\bigr]+\frac{C'(x,x)}{n} +o \biggl(\frac{1}{n}
\biggr).
\end{eqnarray*}
It can be checked, for example, in the case where $W$ is Gaussian that
$C'(x,x)=0$ since moments expand as a series in $1/n^2$. This completes
the proof of the
proposition.
\end{pf}
\subsection{An estimate for $z_c-w_c$ and the role of the fourth moment}
We follow the route developed for Wishart matrices, showing first that
the fluctuations of $z_c^{\pm}(u_0)$ around $w_c^{\pm}(u_0)$ depend on
the fourth moment of the entries of $W$.

%
\begin{proposition}\label{fluctGUEzc}
One has that
\[
\mathbb{E}\bigl[z_c^+(u_0)-w_c^+(u_0)
\bigr]=
\frac{\beta_4m_{\mathrm{sc}}(w_c^+(u_0))^4/(16n)}{(a^{-2}-m_{\mathrm{sc}}'(w_c^+(u_0)))(w_c+m_{\mathrm{sc}}
(w_c^+(u_0))/2)}+o \biggl(\frac{1}{n} \biggr).
\]
As a consequence, for any $\varepsilon>0$ uniformly on $u\in[-2\sigma
+\varepsilon,2\sigma-\varepsilon]$,
%
\begin{equation}
\label{est} \rho_n(u_0)= \sigma_{\mathrm{sc}}^{2\sigma
}(u_0)+
\kappa_4 \frac
{D (u_0)}{n}+o \biggl(\frac{1}{n} \biggr),
\end{equation}
where
$D(u_0)=d(w_c^+(u_0))$ is given for $z\in\mathbb C\setminus\mathbb
R$ by
%
\begin{equation}
\label{defCteu}
d(z)=\frac{1}{16\pi a^2}\Im \biggl( \frac{ m_{\mathrm{sc}}(z)^4}{(a^{-2}-m_{\mathrm{sc}}'(z))(z+m_{\mathrm{sc}}(z)/2)} \biggr).
\end{equation}
%
\end{proposition}

\begin{pf}
We first relate the
critical points $z_c^+$ and $w_c^+$ to the difference of the Stieltjes
transforms $m_n-m_{\mathrm{sc}}$.
The true and approximate critical points satisfy the following equations:
\[
\frac{z_c^+-u_0}{a^2}-m_n\bigl(z_c^+\bigr)=0; \qquad
\frac
{w_c^+-u_0}{a^2}-m_{\mathrm{sc}}\bigl(w_c^+\bigr)=0.
\]
Hence,
%
\begin{equation}
\label{lll} z_c^+-w_c^+=\frac{ (m_n(w_c^+)-m_{\mathrm{sc}}(w_c^+))}{({1}/{a^2}-m_{\mathrm{sc}}'(w_c^+))}+o
\biggl(\frac{1}{n} \biggr),
\end{equation}
where we have used that $m_n-m_{\mathrm{sc}}$ is of order $\frac{1}{n}$. Indeed,
the estimate will again rely on the estimate of the mean of the central
limit theorem for Wigner matrices; see \cite{BaiSil}, Lemma~9.5.
They find that for $z\in\mathbb{C}\setminus\mathbb{R}$,
\[
\lim_{n\to\infty} n\mathbb{E}\bigl(m_n(z)-m_{\mathrm{sc}}(z)
\bigr)= \bigl(1+\tfrac
{1}{4}m_{\mathrm{sc}}'(z)\bigr)
\kappa_4 m_{\mathrm{sc}}(z)^5/16.
\]
Hence, we deduce that
\[
\lim_{n\to\infty} n\bigl(z_c^+-w_c^+
\bigr)=\biggl(\frac{1}{a^2}-m_{\mathrm{sc}'}
\bigl(w_c^+\bigr)\biggr)^{-1} \biggl(1+
\frac{1}{4}m_{\mathrm{sc}}'\bigl(w_c^+\bigr)
\biggr)\kappa_4 m_{\mathrm{sc}}\bigl(w_c^+
\bigr)^5/16.
\]
Using Proposition~\ref{proprhon}, the expansion for the one point
correlation function follows as
\begin{eqnarray*}
\rho_n(x)&=& \sigma^{2\sigma}_{\mathrm{sc}}(u_0)+
\sigma^{2\sigma
}_{\mathrm{sc}}(u_0)\frac
{1}{\Im w_c^+(u_0)}\Im
\bigl(\mathbb E\bigl[z_c^+(u_0)- w_c^+(u_0)
\bigr] \bigr) +o(1/n)
\\
&=&\sigma^{2\sigma}_{\mathrm{sc}}(u_0)+ \frac{\sigma^{2\sigma
}_{\mathrm{sc}}(u_0)}{16\Im
w_c^+(u_0)}
\Im \biggl( \frac{(1+({1}/{4})m_{\mathrm{sc}'}(w_c^+))}{(({1}/{a^2})-m_{\mathrm{sc}'}(w_c^+))} \kappa_4 m_{\mathrm{sc}}
\bigl(w_c^+\bigr)^5 \biggr)
\\
&&{}+o(1/n)
\end{eqnarray*}
It is then an easy computation, using that $m_{\mathrm{sc}}(z)+\frac
{1}{z+m_{\mathrm{sc}}(z)/4}=0$, that this yields Proposition~\ref{fluctGUEzc}.
\end{pf}

\subsection{The localization of eigenvalues}

We now use (\ref{est}) to obtain a precise localization of eigenvalues
in the bulk of the spectrum. A conjecture of Tao and Vu (more precisely
Conjecture~1.7 in \cite{TaoVu}) states that (when the variance of the
entries of $W$ is $\frac{1}{4}$), there exists a constant $c>0$ and a
function $x\mapsto C'(x)$ independent of $\kappa_4$ such that
%
\begin{eqnarray}
\mathbb{E} (\lambda_i- \gamma_i ) &=&
\frac{1}{n\sigma_{\mathrm{sc}}^{2\sigma}(\gamma_i)} \int_{0}^{\gamma
_i}C'(x)\,dx+
\frac{\kappa_4}{2n} \bigl(2\gamma_i^3-
\gamma_i \bigr)+O \biggl(\frac{1}{n^{1+c}} \biggr),
\end{eqnarray}
where $\gamma_i$ is given by
$N_{\mathrm{sc}}(\gamma_i)=i/n$ if $N_{\mathrm{sc}}(x)=\int_{-\infty}^x d\sigma
_{\mathrm{sc}}^{2\sigma}(u)$.
We do not prove the conjecture but another version instead.
More precisely, we obtain the following estimate. Fix $\delta>0$ and
an integer $i$ such that $\delta<i/n<1-\delta$.
Define also
%
\begin{eqnarray}
N_n(x)&:= &\frac{1}{n}\sharp \lbrace i, \lambda_i
\leq x \rbrace\qquad \mbox{with }\lambda_1\leq\lambda_2 \leq
\cdots\leq \lambda _n.
\end{eqnarray}
Let us define the quantile $\hat\gamma_i$ by
\[
\hat\gamma_i:=\inf \biggl\{y, \int^{y}_{-\infty}
\rho_n(x) \,dx=\frac
{i}{n} \biggr\}.
\]
By definition $\mathbb{E} N_n(\hat\gamma_i)=i/n$.
We prove the following result.
%
\begin{proposition}\label{locprop}
Let $\varepsilon>0$ and take $i\in
[\varepsilon n,
(1-\varepsilon)n]$.
There exists a constant $c>0$ such that
\begin{eqnarray}
\hat\gamma_i-\gamma_i &=&\frac{\kappa_4}{2n} \bigl(2
\gamma_i^3-\gamma_i \bigr)+O \biggl(
\frac
{1}{n^{1+c}} \biggr).
\end{eqnarray}
\end{proposition}
The main step to prove this proposition is the following.
%
\begin{proposition}\label{proprole}
Let $\varepsilon>0$. Assume that
$i\in[\frac
{n}{2}, (1-\varepsilon)n] $ without loss of generality. There exists a
constant $c>0$ such that
%
\begin{eqnarray}
&&\hat\gamma_i-\gamma_i -\hat\gamma_{[n/2]}+
\gamma _{[n/2]}
\nonumber
\\[-8pt]
\\[-8pt]
\nonumber
&&\qquad=\frac
{1}{\sigma_{\mathrm{sc}}^{2\sigma}(\gamma_i)}\int_{\gamma_{[n/2]}}^{\gamma
_i}
\bigl[\rho_n(x)-\sigma_{\mathrm{sc}}^{2\sigma}(x) \bigr]\,dx +O
\biggl(\frac
{1}{n^{1+c}} \biggr).
\end{eqnarray}
Note here that $\gamma_{[n/2]}=0$ when $n$ is even.
\end{proposition}

\begin{pf}
 The proof is
divided into Lemmas~\ref{Lemini} and~\ref{Lemreplace} below.
\end{pf}
%
\begin{lemme}\label{Lemini} For any $\varepsilon>0$, there exists
$c>0$ such
that uniformly on $i\in[\varepsilon N,(1-\varepsilon)N]$
%
\begin{equation}
\gamma_i-\hat\gamma_i=\frac{1}{\sigma_{\mathrm{sc}}^{2\sigma}(\gamma_i)} \mathbb{E}
\bigl( N_n(\hat\gamma_i)-N_{\mathrm{sc}}(\hat
\gamma_i) \bigr)+O \biggl(\frac
{1}{n^{1+c}} \biggr).
\end{equation}
\end{lemme}
\begin{pf}
Under assumptions of subexponential tails, it is proved in \cite{EYY2}
and \cite{ErdosRigid} (see also Remark~2.4 of \cite{TaoVu})
that given $\eta>0$ for $n$ large enough
%
\begin{equation}\label{taovuconc}
\mathbb{P} \Bigl( \max_{\varepsilon N\le i\le(1-\varepsilon)n} |\gamma_i-
\lambda_i| \geq n^{\eta-1} \Bigr) \leq n^{-\log n}.
\end{equation}
%
This implies that
\[
\mathbb{E} \bigl[N_n\bigl(\gamma_i+n^{\eta-1}
\bigr) \bigr]\ge\frac
{i}{n}+n^{-\log
n}, \qquad\mathbb{E}
\bigl[N_n\bigl(\gamma_i-n^{\eta-1}\bigr)\bigr]\le
\frac
{i}{n}-n^{-\log
n},
\]
from which it follows that for $n$ large enough
%
\begin{equation}
\label{eq2}
\max_{\varepsilon N\le i\le(1-\varepsilon
)n} |\gamma_i-\hat
\gamma_i| \leq2n^{\eta-1}.
\end{equation}

From the fact that $\mathbb{E}N_n(\hat\gamma_i)= N_{\mathrm{sc}}(\gamma_i)$,
we deduce that
%
\begin{eqnarray}
\mathbb{E}N_n(\hat\gamma_i)-N_{\mathrm{sc}}(\hat
\gamma _i) &=& N_{\mathrm{sc}}(\gamma_i)-N_{\mathrm{sc}}(
\hat\gamma_i)
\nonumber
\\[-8pt]
\\[-8pt]
\nonumber
&=& N_{\mathrm{sc}}'(
\gamma_i) (\gamma_i-\hat\gamma_i ) -\int
_{\gamma
_i}^{\hat\gamma_i}\!\int_{\gamma_i}^u
N_{\mathrm{sc}}''(s) \,ds \,du.
\end{eqnarray}
Using that $N_{\mathrm{sc}}'(x)= \frac{1}{2\pi\sigma^2}\sqrt{4\sigma
^2-x^2}1_{|x|\leq2\sigma}$ and that both $\gamma_i$ and $\hat\gamma
_i$ lie within $(-2\sigma+\varepsilon, 2\sigma- \varepsilon)$ for some
$0<\varepsilon<2\sigma$,
we deduce that
\[
\mathbb{E}N_n(\hat\gamma_i)-N_{\mathrm{sc}}(\hat
\gamma_i)=\sigma _{\mathrm{sc}}^{2\sigma}(\gamma
_i) (\gamma_i-\hat\gamma_i )+ O (
\gamma_i-\hat \gamma _i )^2.
\]

We now make the following replacement.
\end{pf}
%
\begin{lemme}\label{Lemreplace}
Let $\varepsilon>0$. There exist a constant $c>0$ such that uniformly
on $i\in[\varepsilon n, (1-\varepsilon)n]$,
%
\begin{eqnarray}
&&\mathbb{E} \bigl( N_n(\hat\gamma_i)-N_{\mathrm{sc}}(
\hat\gamma_i) \bigr) =\mathbb{E} \bigl(N_n(
\gamma_i)-N_{\mathrm{sc}}(\gamma_i) \bigr)+O \biggl(
\frac
{1}{n^{1+c}} \biggr).
\end{eqnarray}
\end{lemme}

\begin{pf}
We write that
%
\begin{eqnarray}
&&\quad\mathbb{E} \bigl(N_n(\hat\gamma_i)-N_{\mathrm{sc}}(
\hat\gamma_i) \bigr)
\nonumber
\\[-8pt]
\label{enleverlambda}
\\[-8pt]
\nonumber
&&\quad\qquad =\mathbb{E} \bigl(N_n(\gamma_i)-N_{\mathrm{sc}}(
\gamma_i) \bigr)+ \mathbb{E} \bigl(N_n(\hat
\gamma_i)-N_n(\gamma_i)-N_{\mathrm{sc}}(
\hat\gamma _i)+N_{\mathrm{sc}}(\gamma_i) \bigr).
\end{eqnarray}
We show that the second term in (\ref{enleverlambda}) is negligible
with respect to $n^{-1}$. In fact, by \eqref{est} and \eqref{eq2}, for
$\varepsilon>0$, there exists $\delta>0$ such
that for any $i\in[\varepsilon n, (1-\varepsilon)n]$,
%
\begin{eqnarray}
\bigl| \mathbb{E} \bigl(N_n(\hat\gamma_i)-N_n(
\gamma _i)-N_{\mathrm{sc}}(\hat\gamma_i)+N_{\mathrm{sc}}(
\gamma_i) \bigr)\bigr| &\leq& \biggl| \int_{\gamma_i}^{\hat\gamma_i}
\bigl(\rho_n(x)-\sigma _{\mathrm{sc}}^{2\sigma} (x)\bigr) \,dx \biggr|
\nonumber
\\[-8pt]
\\[-8pt]
\nonumber
&\leq & n^{\eta-1}\frac{1}{n}\leq\frac{1}{n^{2-\eta}}.
\end{eqnarray}
In the last line, we have used \eqref{est}.
This completes the proof of Lemma~\ref{Lemreplace}.

Combining Lemmas~\ref{Lemini} and~\ref{Lemreplace} yields
Proposition~\ref{proprole}:
%
\begin{eqnarray}
\delta_n (i)&:=&\sigma_{\mathrm{sc}}^{2\sigma}(
\gamma_i) (\gamma_i-\hat \gamma_i )-
\sigma_{\mathrm{sc}}^{2\sigma}(\gamma_{[n/2]}) (
\gamma_{[n/2]}-\hat\gamma_{[n/2]} )
\nonumber\\
&=& \int
_{\gamma_{[n/2]}}^{\gamma_i} \bigl[\rho_n(x)-\sigma
_{\mathrm{sc}}^{2\sigma}(x) \bigr] \,dx+O \biggl(\frac{1}{n^{1+c}} \biggr)
\nonumber
\\[-8pt]
\label{36}
\\[-8pt]
\nonumber
&=& \frac{1}{n} \int_{\gamma_{[n/2]}}^{\gamma_i}
\kappa_4 D(x) \,dx+O \biggl(\frac{1}{n^{1+c}} \biggr)
\\
&=&
\frac{1}{n} \int_{0}^{\gamma_i}
\kappa_4 D(x) \,dx+O \biggl(\frac
{1}{n^{1+c}} \biggr),\nonumber
\end{eqnarray}
where we used that $\gamma_{[n/2]}$ vanishes or is at most of order $1/n$.
This formula will be the basis for identifying the role of $\kappa_4$
in the $\frac{1}{n}$ expansion of $\hat\gamma_i$.
We now write for a point $x$ in the bulk $(-2\sigma(1-\delta),
2\sigma
(1-\delta))$ that
\[
x= 2\sigma\cos\theta.
\]
We also write that $\gamma_i=2\sigma\cos\theta_0$.
We then have that
\[
w_c(x)=\frac{\cos\theta}{2\sigma}+ \frac{a^2}{\sigma}e^{\pm
i\theta
};\qquad
m_{\mathrm{sc}}\bigl(w_c(x)\bigr)=\pm i \pi\sigma_{\mathrm{sc}}^{2\sigma}(x)-
\frac
{1}{2\sigma^2}x.
\]
By combining Proposition~\ref{fluctGUEzc} and (\ref{defCteu}), we
have that
%
\begin{equation}
D(x)= \Im \biggl( \frac{ m_{\mathrm{sc}}(w_c(x))^4}{16
(w_c(x)+m_{\mathrm{sc}}(w_c(x)))\pi}\bigl(1+o(1)\bigr) \biggr).
\end{equation}

When $a\to0$, we then have the following estimates:
\begin{eqnarray*}
x &\sim & \cos\theta;\qquad  m_{\mathrm{sc}}\bigl(w_c(x)\bigr)
\sim-2e^{-i\theta};\\
 \sigma (x) &\sim & \frac{2}{\pi}\sin\theta;\qquad
w_c+m_{\mathrm{sc}}(w_c)/2\sim i \sin \theta.
\end{eqnarray*}
Using (\ref{36}) and identifying the term depending on $\kappa_4$ in
the limit $a\to0$, we then find that
%
\begin{eqnarray}
\delta_n(i)&=& \frac{1}{n} \int
_{0}^{\gamma_i} \kappa_4 D(x) \,dx+O \biggl(
\frac{1}{n^{1+c}} \biggr)
\nonumber\\
\label{lalala}
&=& \frac{\kappa_4}{n}\int_{\theta_0}^{\pi/2}
\frac{\cos(4 \theta
)}{\pi
}\,d\theta+O \biggl(\frac{1}{n^{1+c}} \biggr)
\\
&=&-
\frac{\kappa_4}{2n}\sigma_{\mathrm{sc}}^{2\sigma}(\gamma_i)\cos
\theta _0 \bigl(2 \cos^2 \theta_0 -1\bigr)+O
\biggl(\frac{1}{n^{1+c}} \biggr),\nonumber
\end{eqnarray}
where in the last line we used that $\frac{1}{4} \sin(4 \theta)=
\sin
\theta\cos\theta\cos(2 \theta)$.
Thus, we have that
%
\begin{eqnarray}\label{eq111}
\delta_n (i)&=&-\frac{\kappa_4}{2n} \sigma_{\mathrm{sc}}^{2\sigma}(
\gamma _i) \bigl(2\gamma_i^3-
\gamma_i\bigr)+O \biggl(\frac{1}{n^{1+c}} \biggr).
\end{eqnarray}

We finally show that
\[
\lim_{n\rightarrow\infty} n ( -\gamma_{[n/2]}+\hat\gamma
_{[n/2]} ) =0
\]
which
completes the proof of Proposition~\ref{locprop}.

To that end, let us first notice that for any
$C^4$ function $f$ whose support is strictly included in
that of $\sigma_{\mathrm{sc}}^{2\sigma}$, we have by \cite{EJP705}, Theorem~1.1,  that
%
\begin{equation}
\label{qwer}
\lim_{n\rightarrow\infty}\mathbb E\Biggl[\sum
_{i=1}^n f(\lambda_i)\Biggr]=m(f)+
\kappa_4 \int_{-1}^1
f(t)T_4(t) \frac
{dt}{\sqrt{1-t^2}}:=m_{\kappa_4}(f),
\end{equation}
with $T_4$ the fourth Tchebychev polynomials and $m(f)$ a linear form
independent of $\kappa_4$.

Next, we can rewrite \eqref{qwer} in terms of the quantiles $\hat
\gamma
_i$ as
\begin{eqnarray*}
m_{\kappa_4}(f)&=& n\int f(x)\rho_n(x)\,dx+o(1)
\\
&= & \sum_{i}f(\hat\gamma_i) +n\sum
_i f'(\hat\gamma_i)\int
_{\hat\gamma
_i}^{\hat\gamma_{i+1}} (x-\hat\gamma_i )
\rho_n(x) \,dx+o(1)
\\
&=& \sum_{i}f(\hat\gamma_i) +
\frac{1}{2} \sum_i f'(\hat
\gamma _i) (\hat\gamma_{i+1}-\hat\gamma_i )
+o(1),
\end{eqnarray*}
where we used that $\hat\gamma_{i+1}-\hat\gamma_i$ is of order
$n^{-1+\eta}$ by \eqref{eq2}. Now, again by \eqref{eq2}, we have
\begin{eqnarray*}
\sum_{i}f(\hat\gamma_i)&=&\sum
_{i}f( \gamma_i)+\sum
_i f'(\gamma _i) (\hat
\gamma_i-\gamma_i )+O \biggl(\frac{1}{n^{-1+2\eta
}} \biggr).
\end{eqnarray*}
Moreover, since $\gamma_{i+1}-\gamma_i$ is at most of order $1/n$, we have
\begin{eqnarray*}
\sum_{i}f( \gamma_i)&=& n\int f(x)
\sigma_{\mathrm{sc}}^{2\sigma}(x)\,dx - \frac
{1}{2}\sum
_i f'( \gamma_i) (
\gamma_{i+1}-\gamma_i ) +o(1).
\end{eqnarray*}
Noting that the first term in the right-hand side vanishes we deduce that
\begin{eqnarray*}
m_{\kappa_4}(f)&=&\sum_{i}f'(
\gamma_i) \biggl[\frac{1}{2}(\hat \gamma _{i+1}-\hat
\gamma_i-\gamma_{i+1}+\gamma_i)+\hat
\gamma_i-\gamma _i \biggr]+o(1),
\end{eqnarray*}
where $\hat\gamma_{i+1}-\hat\gamma_i-\gamma_{i+1}+\gamma_i$ is at most
of order $n^{-1-c}$ by \eqref{eq111} since it is approximately equal
to $(\delta_n(i)-\delta_n(i+1))/\sigma^{2\sigma}_{\mathrm{sc}}(\gamma_i)$.
Hence, we find by \eqref{eq111} that
\begin{eqnarray*}
-m_{\kappa_4}(f) &=&\sum f'(\gamma_i) (
\gamma_i-\hat\gamma_i )+o(1)
\\
&=&\frac{1}{n} \sum_i f'(
\gamma_i) \frac{\sigma_{\mathrm{sc}}^{2\sigma
}(\gamma
_{[n/2]})}{\sigma_{\mathrm{sc}}^{2\sigma}(\gamma_i)} \bigl[n(\gamma _{[n/2]}-\hat
\gamma_{[n/2]}) \bigr]
\\
&&{}+\frac{1}{n}\sum_i \frac{f'(\gamma_i) }{\sigma_{\mathrm{sc}}^{2\sigma
}(\gamma_i)}
\int_0^{\gamma_i} C'(x) \,dx +
\frac{\kappa_4}{2 n}\sum_i f'(
\gamma_i) \bigl(2\gamma_i^3-\gamma
_i \bigr)+o(1)
\\
&=& \sigma_{\mathrm{sc}}^{2\sigma}(\gamma_{[n/2]})\bigl[ n(
\gamma_{[n/2]}-\hat \gamma _{[n/2]})\bigr]\int
_{-2\sigma}^{2\sigma} f'(x) \,dx \\
&&{}+\int
f'(x) \int_0^x
C'(y)\,dy \,dx
\\
&& {}+\frac{\kappa_4 }{2}\int f'(x) \bigl(2x^3-x \bigr)
\sigma _{\mathrm{sc}}^{2\sigma}(x) \,dx+o(1).
\end{eqnarray*}
We finally take $f'$ even, that is $f$ odd in which case the last term
in $\kappa_4$ vanishes,
as well as the term depending on $\kappa_4$ in $m_{\kappa_4}(f)$ as
$T_4$ is even and $f$ odd.
Moreover, $ \sigma_{\mathrm{sc}}^{2\sigma}(\gamma_{[n/2]})$ goes to $1/2$.
Hence, we deduce that there exists a constant independent of $\kappa_4$
such that
\[
\lim_{n\rightarrow\infty} n (\gamma_{[n/2]}-\hat\gamma
_{[n/2]} )=C .
\]
In fact, this constant must vanish as in the case where the
distribution is symmetric, and $n$ even, both
$\gamma_{[n/2]}$ and $\hat\gamma_{[n/2]}$ vanish by symmetry.
\end{pf}

\section*{Acknowledgments}
We gratefully acknowledge the early work with Po-Ru Lu that investigated
real versions of these results. We are indebted to Bernie Wang for so very
much: the careful computations in Julia, incredibly helpful conversations
and his checking the consistency of various formulations.
We thank Peter Forrester for helping us track down results for the
smallest singular value. We also thank Folkmar Bornemann for teaching
us how to use
the beautiful codes in \cite{bornemann2010} for drawing some of these
distributions. We thank an anonymous referee for his careful reading
and helpful remarks.





%

\printaddresses
\end{document}